\documentclass{amsart}

\usepackage[english]{babel}
\usepackage[letterpaper,top=2cm,bottom=2cm,left=3cm,right=3cm,marginparwidth=1.75cm]{geometry}
\usepackage{amssymb}
\usepackage{amsmath}
\usepackage{amsthm}
\usepackage{thmtools}
\usepackage{graphicx}
\IfFileExists{physics.sty}{\usepackage{physics}}{}
\usepackage{tikz-cd}
\usepackage{mathtools}
\usepackage[all]{xy}
\usepackage{xcolor}

\theoremstyle{plain}
\newtheorem{theorem}{Theorem}[section] \newtheorem{corollary}[theorem]{Corollary}
\newtheorem{lemma}[theorem]{Lemma} \newtheorem{proposition}[theorem]{Proposition}
\theoremstyle{definition}
\newtheorem{definition}[theorem]{Definition} \newtheorem{remark}[theorem]{Remark} 

\usepackage[colorlinks=true,allcolors=blue]{hyperref}
\usepackage{cleveref}

\DeclareMathOperator{\Stab}{Stab} \DeclareMathOperator{\gldim}{gldim} \DeclareMathOperator{\Hom}{Hom}
 \DeclareMathOperator{\Aut}{Aut}

\DeclareMathOperator{\ch}{ch} 
 \DeclareMathOperator{\rk}{rk}

\newcommand{\D}{\mathcal{D}} \newcommand{\Db}{\mathrm{D}^{\mathrm b}} \newcommand{\cP}{\mathcal{P}}
\newcommand{\A}{\mathcal{A}}   
  
 \newcommand{\OO}{\mathcal{O}} \newcommand{\C}{\mathbb{C}}
\newcommand{\R}{\mathbb{R}} \newcommand{\Z}{\mathbb{Z}} \newcommand{\QQ}{\mathbb{Q}}
\newcommand{\HH}{\mathbb{H}} \newcommand{\PP}{\mathbb{P}}
 \newcommand{\ii}{\mathrm{i}} 
 \newcommand{\Sfun}{\mathbb{S}}

\newcommand{\Q}{\QQ} \newcommand{\Pp}{\PP} \newcommand{\Perf}{\operatorname{Perf}}
\newcommand{\Ind}{\operatorname{Ind}} \newcommand{\Ku}{\operatorname{Ku}}
\newcommand{\Knum}{\mathrm K_{\mathrm{num}}} \newcommand{\MF}{\operatorname{MF}}
\newcommand{\thick}{\operatorname{thick}} \newcommand{\RHom}{\operatorname{RHom}}
 \newcommand{\id}{\operatorname{id}} \newcommand{\gr}{\mathrm{gr}}
\newcommand{\sg}{\mathrm{sg}} 
\hypersetup{pdftitle={Construction of stability conditions on some Kuznetsov
components},pdfauthor={Dongjian Wu}}
\title{Construction of stability conditions on some Kuznetsov components}
\author[D. Wu]{Dongjian Wu}
\address{D. Wu: Shanghai Institute for Mathematics and
Interdisciplinary Sciences (SIMIS), Shanghai 200433, China \& Research Institute
of Intelligent Complex Systems, Fudan University, Shanghai 200433, China}
\email{wdj@simis.cn}

\date{}
\begin{document}
\begin{abstract}
We construct Bridgeland stability conditions on the Kuznetsov
components of smooth Fano hypersurfaces of degrees three, four,
and six and dimension at least three, smooth quartic double
fivefolds, and smooth Fano complete intersections of three
quadrics. For odd-dimensional cubic hypersurfaces and quartic
double fivefolds, every connected component of the numerical
stability manifold is biholomorphic to $\C\times\HH$. For general cubic
sevenfolds and quartic double fivefolds, spherical twists
produce infinitely many distinct connected components.

\end{abstract}
\maketitle
\section{Introduction}

Bridgeland stability conditions generalize slope stability for
coherent sheaves to triangulated categories, and their parameter
spaces carry a natural complex structure \cite{Bridgeland}.
They provide a way to study stable objects and their moduli spaces
in categories that need not arise as derived categories of varieties.
An important class of examples consists of the Kuznetsov components
of Fano varieties. In this paper, we construct stability conditions
on several families of these components and describe their numerical
stability manifolds in certain cases.

\subsection{Background and motivations}

For a Fano variety, a semiorthogonal decomposition often separates
a collection of exceptional objects from a residual category.
The latter is called a Kuznetsov component. Early studies of these
categories revealed equivalences between the residual categories
of cubic and $V_{14}$ threefolds \cite{KuzCubicThreefold2004}.
Their systematic study for Fano threefolds and cubic fourfolds
was developed in \cite{KuzFanoThreefolds2009,KuzCubicFourfolds2010}.
These works, together with \cite{IKP2014}, connect residual
categories with classical questions about Fano varieties, including
rationality and the relation between different families.
For cubic fourfolds, the conjecture of
\cite{KuzCubicFourfolds2010} relates rationality to whether the
Kuznetsov component is equivalent to the derived category of a K3
surface.

Kuznetsov components often have Calabi--Yau or fractional
Calabi--Yau structure of smaller dimension than the ambient
variety. This motivates the study of stability conditions on
them, where moduli spaces and wall-crossing may be easier to
control. Moduli spaces of stable objects can recover classical
parameter spaces, produce new varieties, and help distinguish
derived categories. This approach has led to concrete results. For cubic threefolds,
stability conditions recover the Fano surface of lines and give a
categorical Torelli theorem \cite{BMMS2012}. Stability conditions
on the Kuznetsov components of the Picard-rank-one Fano threefolds
covered by \cite{BLMS2023} and of smooth cubic fourfolds were
constructed there using a criterion for restriction to the
orthogonal complement of an exceptional collection. The
construction extends to families and to Gushel--Mukai varieties
\cite{BLMNPS2021,PerryPertusiZhao2022}. In the K3 cases, moduli
spaces of stable objects recover classical hyperk\"ahler varieties
and produce new ones \cite{LPZ2023,BLMNPS2021}. Other applications
include moduli spaces of sheaves, conics and instantons
\cite{LiuZhang2022,Zhang2021}, as well as Lagrangian constructions
obtained by restriction to hyperplane sections
\cite{FGLZ2025,LLPZ2026}.

The cubic fivefold extends this picture to fractional
Calabi--Yau dimension $7/3$. For every smooth cubic fivefold,
\cite[Theorem~1.1]{PeizeLiu2025} constructs a family of
Serre-invariant stability conditions in a single
$\widetilde{\mathrm{GL}}_2^+(\R)$-orbit containing a Gepner
condition for the rotation functor. The construction uses a
quadric surface fibration, a Bogomolov--Gieseker inequality for
Clifford modules on $\PP^3$, and the restriction criterion of
\cite{BLMS2023}. For a general cubic fivefold, the same work proves
non-emptiness of moduli spaces for every nonzero numerical class
and gives Lagrangian restriction constructions for primitive
classes \cite[Theorem~1.4]{PeizeLiu2025}. These results motivate
the existence problem beyond the threefold and K3 cases.

Existence on a Kuznetsov component requires additional work even
when stability conditions are available on the ambient variety.
The construction of \cite{Li2026} gives stability conditions on
$\Db(X)$, and \cite{KLP2026} induces bounded $t$-structures on
semiorthogonal components. A stability condition on a prescribed
component still requires a compatible central charge, the
Harder--Narasimhan property and the support property. For cubic
sevenfolds, \cite[Theorem~1.1]{PeizeLiu2026} embeds the Kuznetsov
component into a Clifford category on $\PP^4$. The proposed
restriction construction requires a further Chern inequality, and
existence is formulated in \cite[Conjecture~1.2]{PeizeLiu2026}.
This suggests looking for a construction that applies uniformly
in the dimension.

Once existence is established, one can ask how much of the stability
manifold the construction describes. A distinguished
Serre-invariant component is contractible for cubic threefolds,
quartic double solids and Gushel--Mukai threefolds
\cite[Theorem~1.1]{FanLiuMa2023}. For very general cubic fourfolds,
the entire numerical stability manifold is a single
$\widetilde{\mathrm{GL}}_2^+(\R)$-orbit
\cite[Corollary~3.9]{FeyzbakhshPertusi2023}. For cubic fivefolds,
uniqueness of the Serre-invariant orbit is raised in
\cite[Remark~5.6]{PeizeLiu2025}. These questions also lead to the
study of connected components and the action of autoequivalences.

\subsection{Main results}

Throughout the paper,
numerical stability means stability with support on the full
numerical Grothendieck group, unless a different lattice is
specified. For a smooth Fano hypersurface
$X_d\subset\PP^{n+1}$ of degree $d$, we use the convention
\[
 \Ku(X_d)=\langle\OO_{X_d},\OO_{X_d}(1),\ldots,
                         \OO_{X_d}(n+1-d)\rangle^\perp.
\]
We denote its rotation functor by $R_X$. Under
\eqref{eq:graded-hypersurface-equivalence}, it corresponds to the
grading shift by \eqref{eq:grading-rotation-compatibility}.
In particular, $R_X^d\simeq[2]$.

Our first result gives stability conditions on three families
of Fano hypersurfaces.

\begin{theorem}
\label{thm:intro-hypersurfaces}
Let $d\in\{3,4,6\}$ and let $X_d\subset\PP^{n+1}$ be a smooth
Fano hypersurface with $n\geq\max\{3,d-1\}$. There exists a numerical
Bridgeland stability condition $\widetilde\sigma_X$ on $\Ku(X_d)$
such that
\[
  R_X\widetilde\sigma_X =\widetilde\sigma_X\cdot\frac{2}{d}.
\]
It has mass-Hom bounds and support on $\Knum(\Ku(X_d))$. Moreover,
\[
 \gldim(\widetilde\sigma_X)=\frac{(n+2)(d-2)}{d}.
\]
\end{theorem}

The cubic case applies in every dimension at least three. It
proves \cite[Conjecture~1.2]{PeizeLiu2026} for cubic sevenfolds
and recovers the fivefold existence and Gepner symmetry of
\cite{PeizeLiu2025}. For quartics and sextics, the same
construction applies in every Fano dimension covered by the
theorem. The rotation equation is the Gepner symmetry considered
in \cite{Toda2013,Toda2016}. 

We also prove existence for complete intersections of three
quadrics.

\begin{theorem}
\label{thm:intro-three-quadrics}
Let $X\subset\PP^{n+3}$ be a smooth Fano complete intersection of three
quadrics. Then
\[
\Stab_{\mathrm{num}}(\Ku(X))\neq\varnothing, \qquad
\Ku(X)=\langle\OO_X,\OO_X(1),\ldots,
   \OO_X(n-3)\rangle^\perp.
\]
The constructed stability condition has support on the full numerical Grothendieck group.
\end{theorem}

For quartic double fivefolds, we obtain existence by a weighted
version of the hypersurface construction in
Proposition~\ref{prop:quartic-double-existence}.

Our next result describes every connected component for
odd-dimensional cubic hypersurfaces and quartic double fivefolds.

\begin{theorem}
\label{thm:intro-cubic-components}
Let $X\subset\PP^{n+1}$ be a smooth cubic hypersurface of odd dimension $n\geq3$. There are
biholomorphisms
\[
 \Stab_{\mathrm{num}}(\Ku(X)) \simeq \coprod_{\alpha\in I_X}(\C\times\HH), \qquad\text{ where } I_X=\pi_0\Stab_{\mathrm{num}}(\Ku(X)).
\]
Every connected component is a free $\widetilde{\mathrm{GL}}_2^+(\R)$-orbit and is
contractible. If $X$ is a general cubic sevenfold, then $I_X$ is infinite.
The same description holds for a smooth quartic double fivefold $Y\subset\PP(1,1,1,1,1,1,2)$, and $I_Y$ is infinite when $Y$ is general. In both
general cases, a spherical twist $T_S$ satisfies
\(
 T_S^m(\Sigma)\neq\Sigma\) {for every component $\Sigma$ and every $m\neq0$.}
\end{theorem}

The infinitude statement applies to general members, where the
spherical bundles of \cite{IlievManivel2014,Abuaf2020} exist.

\subsection{Strategy of proof}

For hypersurfaces, we work with adjunctions between graded
singularity categories and then transport stability conditions to
Kuznetsov components through
\eqref{eq:factorization-singularity-equivalence} and
\eqref{eq:graded-hypersurface-equivalence}.

For a coordinate of weight $a$ in a degree-$d$ potential,
the relevant exact triangle expresses $VU(E)$ as an extension
of objects of phases $\phi$ and $\phi+1-2a/d$, where $E$ has
phase $\phi$. Thus $a\leq d/2$ gives the phase inequality
required by \cite[Appendix~A]{LLL2026} to construct
$t$-structures on the target category. Conservativity and
reflection of compact objects ensure that these $t$-structures
are bounded and give finite Harder--Narasimhan filtrations. 

The transfer of support requires a generator argument. The
derived self-intersection triangle may have a nonzero connecting
morphism, but properness and graded periodicity imply that a
power of this morphism vanishes. Repeated use of the octahedral
axiom then shows that derived extension carries a classical
generator to a classical generator. This transfers mass-Hom
bounds in the sense of \cite[Definition~2.7]{HLR2025}, extending
the adjunction estimates of \cite[Lemma~7.4]{LLL2026}.
An Euler-functional argument, following
\cite[Theorem~2.1]{Cheng2026}, gives support on the full numerical
group without requiring an injective numerical pushforward.

The starting stability conditions come from the elliptic diagonal types
$(3,3,3)$, $(4,4,2)$ and $(6,3,2)$. Product constructions
\cite{LLL2026,LMP2025,Liu2021}, finite character descent
\cite{BFK2014}, and a reverse finite-cover mass estimate give
Gepner conditions in arbitrarily large dimensions. Coordinate
restrictions remove the weighted variables in degrees four and
six. A decomposition into powers of linear forms realizes any
smooth equation as a linear section of a resulting Fermat
hypersurface. A flag of smooth sections then supplies the
required sequence of restrictions. The same method, retaining
one variable of weight two, treats quartic double fivefolds.

For complete intersections of three quadrics, the construction
depends on the parity of the dimension. In even dimension, we
use the twisted-surface description of
\cite[Proposition~8.7(ii)]{CO2026}. We show that its geometric
hypotheses follow from smoothness of the complete intersection
and construct a twisted surface stability condition with full
numerical support. In odd dimension, we use the Clifford category
on $\PP^2$ from \cite[Theorem~5.5]{KuzQuadrics2008}. A weak
discriminant estimate gives the inequality needed for a surface
tilt. The Chern character of the underlying ordinary sheaf
identifies the full rational numerical group with the space of
rank, degree and second Chern character coordinates.

For the component description, we use numerical rank two, a
classical generator, and the failure of finite generation of
$K_0$. A central charge of real rank one would give a
finite-length heart. The generator would force this heart to
have finitely many simple objects, contradicting the failure of
finite generation of $K_0$. Every central charge therefore has
real rank two. The
$\widetilde{\mathrm{GL}}_2^+(\R)$-orbits are then open and are
the connected components. The central-charge map identifies
each component with the universal cover of $\C^*\times\HH$,
which gives the biholomorphism with $\C\times\HH$.

Finally, if a positive power of a spherical twist preserved a
component, its phase transformation would decrease every phase.
The twist triangles would then place every semistable object in
the subcategory generated by the spherical object. Finite
Harder--Narasimhan filtrations would imply that this subcategory
is the whole category. The positive dg-algebra theorem of
\cite[Theorem~8.1]{KellerNicolas2013} rules this out in numerical
rank two. Taking inverses gives the assertion for negative
powers.

\subsection*{Contents}
In Section~\ref{sec:preliminaries}, we recall stability
conditions, the global dimension function and graded singularity
categories. In Section~\ref{sec:coordinate-induction}, we study
coordinate restrictions between weighted hypersurface residual
categories. In Section~\ref{sec:elliptic-seeds}, we construct
 stability conditions on kuznetsov components of Fano hypersurfaces of degrees three, four, and six from products of elliptic
curves. In Section~\ref{sec:hypersurface-applications}, we apply
these constructions to ordinary hypersurfaces and quartic
double fivefolds. In Section~\ref{sec:three-quadrics}, we treat
complete intersections of three quadrics. In
Section~\ref{sec:stability-components}, we describe connected
components and study the action of spherical twists.

\subsection*{Acknowledgements}
The author thanks Yu-Wei Fan for his continuous support,
and Nantao Zhang and Tianle Mao for many helpful discussions.
\subsection*{AI disclosure}
We acknowledge the use of AI in preparing this manuscript.
With the help of ChatGPT, the author initially explored
the use of a suitable function to develop a categorical
construction of stability conditions on
noncommutative varieties, motivated by the induction
construction in \cite{Li2026}. This approach has not yet
been successful. During this investigation, however, the
author found a different method for constructing stability
conditions on the Kuznetsov components considered in this
paper. ChatGPT helped accelerate this work by assisting
in filling gaps in some proofs, checking arguments, and
editing the language. The author takes full responsibility for all mathematical
claims, arguments, and conclusions, as well as for any
remaining errors.

\section{Preliminaries}
\label{sec:preliminaries}

In this section, we recall Bridgeland stability conditions and mass-Hom bounds
\cite{Bridgeland,Bayer2019,HLR2025}, the global dimension function \cite{IQ2023,Qiu},
and the Kuznetsov components and graded factorizations used below
\cite{Orlov2009,KP2021,Toda2013}.

\subsection{Bridgeland stability conditions}
Let $\D$ be a triangulated category and let $K_0(\D)$ be its Grothendieck group. Fix a free abelian group $\Lambda$ of finite
rank and a surjective class map
\(
v:K_0(\D)\twoheadrightarrow\Lambda.
\)

\begin{definition}[{\cite[Definition 3.3]{Bridgeland}}]
A \emph{slicing} $\cP$ on $\D$ is a collection of full additive
subcategories $\cP(\phi)\subset\D$, indexed by $\phi\in\R$,
satisfying:
\begin{enumerate}
\item $\cP(\phi+1)=\cP(\phi)[1]$ for every $\phi\in\R$;
\item if $\phi_1>\phi_2$, $E_i\in\cP(\phi_i)$, then
$\Hom_{\D}(E_1,E_2)=0$;
\item every nonzero object $E\in\D$ admits a finite sequence of distinguished
triangles \[0 = \xymatrix @C=5mm{ E_0 \ar[rr] && E_1 \ar[dl] \ar[rr] && E_2 \ar[dl] \ar[r] & \dots \ar[r] & E_{m-1} \ar[rr]
&& E_m \ar[dl] \\
& A_1 \ar@{-->}[ul] && A_2 \ar@{-->}[ul] &&&& A_m \ar@{-->}[ul] } = E
\]
 such that each nonzero factor
    $A_i\in\cP(\phi_i)$ with real numbers $\phi_1>\phi_2>\cdots>\phi_m$.
\end{enumerate}
The sequence in the last condition is called the \emph{Harder--Narasimhan filtration} (\emph{HN-filtration}) of $E$, and the $A_i$ are its HN
factors.
\end{definition}

For $0\neq E\in\D$, we denote the largest and smallest phases of its HN factors by
\[
  \phi^+_{\cP}(E)=\phi_1, \qquad \phi^-_{\cP}(E)=\phi_m.
\]
For an interval $I\subset\R$, let $\cP(I)$ be the extension-closed subcategory generated by $\cP(\phi)$ for $\phi\in I$. In particular,
$\A_{\cP}:=\cP((0,1])$ is the heart of a bounded $t$-structure. Given two slicings $\cP_1$ and $\cP_2$ on $\D$, their distance is defined as
\[
d(\cP_1,\cP_2)
:=
\sup_{0\neq E\in\D}
\left\{
\left|\phi^+_{\cP_1}(E)-\phi^+_{\cP_2}(E)\right|,
\left|\phi^-_{\cP_1}(E)-\phi^-_{\cP_2}(E)\right|
\right\}.
\]

We use the norm formulation of the support property from \cite{KS2008}. Its relation with the quadratic
formulation and the corresponding deformation theorem are discussed in \cite{Bayer2019}.

\begin{definition}[{\cite[Definition 1.1]{Bridgeland}}]
A \emph{Bridgeland pre-stability condition} on $\D$ with respect to
$(\Lambda,v)$ is a pair $\sigma=(\cP,Z)$ consisting of a group homomorphism
$Z:\Lambda\to\C$ and a slicing $\cP$ such that, for every
$0\neq E\in\cP(\phi)$,
\[
Z(v(E))=m(E)e^{\pi\ii\phi} \qquad\text{for some }m(E)>0.
\]
A nonzero object $E$ is called $\sigma$\emph{-semistable} of phase $\phi$ if \(E\in\cP(\phi)\). It is called $\sigma$\emph{-stable} if it is simple in $\cP(\phi)$. A pre-stability condition $\sigma=(\cP,Z)$ is called locally finite if for any $\phi$, there exists an
$\epsilon>0$ such that the quasi-abelian category $\cP(\phi-\epsilon,\phi+\epsilon)$ is finite length. We say a pre-stability condition $\sigma=(\cP,Z)$ satisfies the \emph{support property} if, for any
norm $\lVert\cdot\rVert$ on $\Lambda_{\R}$, there exists a constant $C>0$ such that
\[
\lVert v(E)\rVert\leq C\lvert Z(v(E))\rvert
\]
for every semistable object $E$. A pre-stability condition satisfying the
support property is called a \emph{stability condition}.
\end{definition}
Denote the set of all stability conditions with respect to $(\Lambda,v)$ by $\Stab_\Lambda(\D)$. When the support lattice is clear from the
context, we simply write $\Stab(\D)$. We also abbreviate $Z(v(E))$ to $Z(E)$ and write $\phi^\pm_\sigma(E):=\phi^\pm_{\cP}(E)$ and
$\A_\sigma:=\cP((0,1])$. Equivalently, a stability condition can be described by a heart of a bounded $t$-structure with a stability function satisfying the
HN property and the support property \cite[Proposition~5.3]{Bridgeland}.

For a smooth projective variety $X$, we write $K_{\rm num}(X)$ for the numerical Grothendieck group of $\Db(X)$. Throughout
this paper, stability conditions on $\Db(X)$ are considered with respect to a finite-rank free abelian group
$\Lambda$ and a class map
\(
v\colon K_0(X)\longrightarrow\Lambda
\)
which factors through $K_{\rm num}(X)$. We call these
\emph{numerical stability conditions}. The lattice $\Lambda$
is either $K_{\rm num}(X)$ or a numerical quotient
specified by the construction. We write $\Stab_\Lambda(X)$
when the support lattice needs to be specified, and
$\Stab(X)$ when it is understood.

For an admissible subcategory $\A\subset\Db(X)$, we similarly set
\[
 \Knum(\A)
 =
 \bigl(K_0(\A)/\operatorname{rad}(\chi_\A)\bigr)
 /\mathrm{torsion},
\]
where $\chi_\A$ denotes the Euler pairing. Stability conditions with
respect to a finite-rank support lattice whose class map factors through
$\Knum(\A)$ are also called numerical. We write $\Stab_{\mathrm{num}}(\A)$ when the
support lattice is the full group $\Knum(\A)$.

For an object $0\neq E$ with HN factors $A_1,\ldots,A_m$, its \emph{mass} is
\[
  m_\sigma(E)=\sum_{j=1}^m|Z(A_j)|.
\]

For stability conditions $\sigma_1=(\cP_1,Z_1)$ and $\sigma_2=(\cP_2,Z_2)$, their distance is defined by
\[
d(\sigma_1,\sigma_2) := \sup_{0\ne E\in\D} \left\{
|\phi^-_{\sigma_2}(E)-\phi^-_{\sigma_1}(E)|,
  |\phi^+_{\sigma_2}(E)-\phi^+_{\sigma_1}(E)|,
\left|\log\frac{m_{\sigma_2}(E)}{m_{\sigma_1}(E)}\right| \right\}.
\]
\begin{theorem}[{\cite[Theorem 1.2]{Bridgeland}}]
\label{thm:deformation}
The space of stability conditions $\Stab(\D)$ has the structure of a complex manifold, and the forgetful map
\[
\Stab(\D)\to\Hom_{\Z}(\Lambda,\C)
\]
that sends a stability condition to its central charge is a local isomorphism.
\end{theorem}

We now recall two standard group actions on $\Stab(\D)$. Let $\Phi$ be an exact autoequivalence of $\D$ preserving
$\ker v$, and denote by $\Phi_{\ast}\colon K_0(\D)\to K_0(\D)$ the induced isomorphism on the Grothendieck group. For a stability condition $\sigma=(\cP,Z)$ on $\D$, we define the left action of $\Phi$ on $\sigma$ by 
\[
  \Phi\cdot\sigma\coloneq(\Phi\circ\cP,Z\circ\Phi_{\ast}^{-1}).
\]
Let $\widetilde{\mathrm{GL}}_2^+(\R)$ be the universal cover of $\mathrm{GL}_2^+(\R)$. For any element $g=(T,f)\in\widetilde{\mathrm{GL}}_2^+(\R)$, the right $\widetilde{\mathrm{GL}}_2^+(\R)$-action is defined by
\[
 \sigma[g]=(\cP[g],Z[g]), \quad Z[g]=T^{-1}Z,\quad \cP[g](\phi)=\cP(f(\phi)).
\]
In particular, for $z=a+\ii b\in\C$, the right $\C$-action on $\sigma$ is defined by
\[
(\cP,Z)\cdot z = \bigl(\cP_z,e^{-\pi\ii z}Z\bigr), \qquad
   \cP_z(\phi)=\cP(\phi+a).
\]

\begin{definition}
Let $R$ be an exact autoequivalence of $\D$ satisfying
$R^d\simeq[2]$. A stability condition
$\widetilde\sigma=(\mathcal P,Z)$ is of \emph{Gepner
type} for $R$ if
\begin{equation}
\label{eq:Gepner}
  R\widetilde\sigma =\widetilde\sigma\cdot\frac{2}{d}.
\end{equation}
\end{definition}

With our conventions for the group actions, \eqref{eq:Gepner} is equivalent to
\[
R\bigl(\mathcal P(\phi)\bigr) =\mathcal P\left(\phi+\frac2d\right), \qquad Z(RE)=e^{2\pi\ii/d}Z(E).
\]
Thus $R$ preserves semistability and increases phases by
$2/d$. Throughout this paper, Gepner type refers to this
symmetry. 

Finally, we recall the notion of mass-Hom bounds for Bridgeland
stability conditions, following \cite{HLR2025}. For the following
definition and lemma, let $\D$ be an idempotent-complete
$\C$-linear triangulated category with a smooth proper dg enhancement.

\begin{definition}[{\cite[Definition~2.7]{HLR2025}}]
\label{def:mass-hom}
A stability condition $\sigma$ on $\D$ has \emph{mass-Hom bounds}
if, for every $A\in\D$, there is a constant $c_A>0$ such that, for all $E\in\D$,
\[
\dim_\C\Hom(A,E)\leq c_A m_\sigma(E).
\]
\end{definition}

We will use the equivalent formulation with total morphism
spaces \cite[Proposition~2.13]{HLR2025}.

\begin{lemma}
\label{lem-total-mass-hom}
A stability condition $\sigma$ on $\D$ has mass-Hom bounds if and only if, for every $A\in\D$, there is a constant $C_A>0$ such that, for all
$E\in\D$, we have
\begin{equation}
 \sum_{j\in\Z}\dim_\C\Hom(A,E[j])
 \leq C_A m_\sigma(E)
 \label{eq-mass-hom-input}
\end{equation}
\end{lemma}

\begin{proof}
It is enough to prove the forward implication. Fix $A\in\D$ and let $\Sfun$ be the Serre
functor. Suppose first that $F$ is semistable. Shifting $F$, we may assume that its phase satisfies
$0<\phi\leq1$. If $\Hom(A,F[j])\neq0$, by Serre duality, we have
\[
  \phi^-(A)\leq\phi+j\leq\phi^+(\Sfun A).
\]
Let $N_A$ be the number of integers in $[\phi^-(A)-1,\phi^+(\Sfun A)]$. Applying the mass-Hom bound to each $F[j]$, we obtain
\[
\sum_{j\in\Z}\dim_\C\Hom(A,F[j]) \leq N_Ac_A m_\sigma(F).
\]

For an arbitrary object $E$, let $F_1,\ldots,F_r$ be its
HN factors. Then, we have
\[
 \begin{aligned}
\sum_{j\in\Z}\dim_\C\Hom(A,E[j]) &\leq \sum_{i=1}^r\sum_{j\in\Z}\dim_\C\Hom(A,F_i[j])\\
&\leq N_Ac_A\sum_{i=1}^r m_\sigma(F_i) =N_Ac_A m_\sigma(E),
 \end{aligned}
\]
which completes the proof.
\end{proof}

\begin{remark}
The property of mass-Hom bounds can be checked on a classical generating set
and on stable objects \cite[Lemma~2.8]{HLR2025}. In addition, it is
preserved at finite Bridgeland distance \cite[Corollary~2.14]{HLR2025}, and
therefore throughout a connected component \cite[Remark~7.3]{LLL2026}.
\end{remark}

\subsection{The global dimension function}
Let $\cP$ be a slicing on a triangulated category $\D$. The
\emph{global dimension} of $\cP$ is defined by
\[
  \gldim(\cP)
  :=
  \sup\bigl\{
    \phi_2-\phi_1
    \bigm|
    \Hom_{\D}\bigl(\cP(\phi_1),\cP(\phi_2)\bigr)\neq0
  \bigr\}
  \in[0,+\infty].
\]
For a stability condition $\sigma=(\cP,Z)$ on $\D$, set $\gldim(\sigma):=\gldim(\cP)$. This is the global dimension function introduced in
\cite{IQ2023,Qiu}. By definition, \(\gldim(\sigma)\) measures the largest phase gap across which a nonzero morphism between semistable objects can occur. The
function $\gldim:\Stab_\Lambda(\D)\to[0,+\infty]$ is continuous and is invariant under exact autoequivalences and the $\C$-action
\cite[Section~3.1]{Qiu}. Hence it descends to a continuous function on the reduced quotient:
\[
\gldim\colon \Aut(\D)\backslash\Stab(\D)/\C\to[0,\infty].
\]

\begin{lemma}[{\cite[Proposition~3.10]{Qiu}}]
\label{lem:serre-phase-gldim}
Suppose that the Serre functor $\Sfun$ takes every $\sigma$-semistable object of phase $\phi$ to a semistable object of phase
$\phi+s$. Then $\gldim(\sigma)=s$.
\end{lemma}

\subsection{Kuznetsov components and graded factorizations}

Let $\D$ be a $\C$-linear triangulated category. A full triangulated subcategory $\A\subset\D$ is \emph{admissible} if its inclusion admits both
left and right adjoints. A \emph{semiorthogonal decomposition}
\[
  \D=\langle\A_1,\ldots,\A_r\rangle
\]
consists of full triangulated subcategories which generate $\D$ and satisfy
\[
\Hom(A_j,A_i[k])=0 \qquad \text{for }j>i,\quad A_j\in\A_j,\quad A_i\in\A_i,
   \quad k\in\Z.
\]
An object $E$ is \emph{exceptional} if
\[
 \Hom(E,E[k])=
 \begin{cases}
\C,&k=0,\\
  0,&k\neq0.
 \end{cases}
\]
A sequence of exceptional objects $(E_1,\ldots,E_r)$ is an exceptional
collection if $\Hom(E_j,E_i[k])=0$ for $j>i$ and all $k\in\Z$.

Let $X$ be a smooth Fano variety, and choose an ample line bundle $H$ and an integer $I_X>0$ such that
\(
\omega_X^{-1}\simeq H^{\otimes I_X}\). Write $\OO_X(t)=H^{\otimes t}$. By Kodaira vanishing theorem, the line bundles
\[
 \OO_X,\OO_X(1),\ldots,\OO_X(I_X-1)
\]
form an exceptional collection.

\begin{definition}
\label{def:kuznetsov-component}
The \emph{Kuznetsov component} associated to $(X,H)$ is defined by
\[
 \begin{aligned}
 \Ku_H(X)
 &\coloneq
 \langle\OO_X,\OO_X(1),\ldots,\OO_X(I_X-1)\rangle^\perp\\
 &=\bigl\{
 E\in\Db(X)\,\bigm|\,
 \Hom(\OO_X(t),E[j])=0
 \text{ for }0\leq t<I_X,\ j\in\Z
 \bigr\}.
 \end{aligned}
\]
We write $\Ku(X)$ when the polarization $H$ is understood.
\end{definition}

There is a semiorthogonal decomposition
\[
 \Db(X)=
 \langle\Ku(X),\OO_X,\OO_X(1),\ldots,\OO_X(I_X-1)\rangle.
\]
In particular, $\Ku(X)$ is admissible and inherits a
smooth proper dg enhancement from $\Db(X)$. For a smooth Fano complete intersection
$X\subset\PP^{n+c}$ of dimension $n\geq1$ and multidegree
$(d_1,\ldots,d_c)$, we take $H=\OO_X(1)$. Adjunction gives
\[
 I_X=n+c+1-\sum_{i=1}^c d_i.
\]
In particular, for a hypersurface
$X_d\subset\PP^{n+1}$ of degree $d$, it gives
\[
 \Ku(X_d)= \langle\OO_{X_d},\OO_{X_d}(1),\ldots,
  \OO_{X_d}(n+1-d)\rangle^\perp,
\]
and for a smooth Fano complete intersection of three quadrics
$X\subset\PP^{n+3}$, it gives
\[
\Ku(X)= \langle\OO_X,\OO_X(1),\ldots,\OO_X(n-3)\rangle^\perp.
\]

For a Fano hypersurface, the \emph{rotation functor} is given by
\[
R_X\coloneq \left. \bigl(\mathbf L_{\OO_X}\circ (-\otimes\OO_X(1))\bigr) \right|_{\Ku(X)},
\]
where $\mathbf L_{\OO_X}$ denotes left mutation through $\OO_X$, defined by the evaluation triangle
\[
  \RHom(\OO_X,F)\otimes_\C\OO_X \longrightarrow F \longrightarrow\mathbf L_{\OO_X}(F).
\]
By \cite[Proposition~3.7]{KP2021}, the composition $\mathbf
L_{\OO_X}\circ(-\otimes\OO_X(1))$ preserves $\Ku(X)$ and restricts to an autoequivalence
$R_X$ of $\Ku(X)$.

Let $P=\C[x_0,\ldots,x_m]$ with $\deg(x_i)=w_i>0$ and let $W\in P$ be a nonzero
homogeneous polynomial of degree $d>0$. For a graded $P$-module $M$, the
\emph{internal grading shift} is defined by
\(
 M(k)_j=M_{j+k}\). We recall the category of graded matrix factorizations,
following \cite[Section~3.1]{Orlov2009} and
\cite[Section~2.2]{Toda2016}.

\begin{definition}
\label{def:graded-matrix-factorization}
A \emph{graded matrix factorization} of $W$ consists of finite-rank graded free $P$-modules $E_0,E_1$ and
$P$-linear maps of degree zero
\[
E_0\xrightarrow{\delta_0}E_1 \xrightarrow{\delta_1}E_0(d)
\]
satisfying
\[
\delta_1\delta_0=W\id_{E_0}, \qquad \delta_0(d)\delta_1=W\id_{E_1}.
\]
Here multiplication by $W$ is regarded as a degree-zero map $E_i\to E_i(d)$.
\end{definition}

For two factorizations $E=(E_0,E_1,\delta_0,\delta_1)$ and
$F=(F_0,F_1,\epsilon_0,\epsilon_1)$, a morphism $f\colon E\to F$ is a pair of
degree-zero maps $f_i\colon E_i\to F_i$ satisfying
\[
  f_1\delta_0=\epsilon_0f_0, \qquad f_0(d)\delta_1=\epsilon_1f_1.
\]
It is null-homotopic if there are degree-zero maps
\[
h_0\colon E_0\to F_1(-d), \qquad h_1\colon E_1\to F_0
\]
such that
\[
 f_0=\epsilon_1(-d)h_0+h_1\delta_0,
 \qquad
 f_1=\epsilon_0h_1+h_0(d)\delta_1.
\]
Denote by $\MF^\Z(W)$ the resulting homotopy
category, which is triangulated, with distinguished triangles
given by mapping cones \cite[Proposition~3.4]{Orlov2009}.
Its translation functor is given by
\[
 E[1]=
 \bigl(
 E_1\xrightarrow{-\delta_1}E_0(d)
 \xrightarrow{-\delta_0(d)}E_1(d)
 \bigr).
\]
The internal grading shift acts on both modules and maps and induces an autoequivalence $(1)$ of $\MF^\Z(W)$. Applying the
translation functor twice gives $E[2]=E(d)$, and hence \((d)\simeq[2]\).

Set $A=P/(W)$. The graded singularity category of $A$ is
\[
  \mathrm D_{\sg}^{\gr}(A) = \Db(\operatorname{grmod}A)/\Perf^{\gr}(A),
\]
where $\operatorname{grmod}A$ denotes finitely generated graded $A$-modules and $\Perf^{\gr}(A)$ denotes graded perfect complexes. Since $W$
annihilates $\operatorname{coker}(\delta_0)$, this cokernel is a graded $A$-module. The cokernel functor induces an equivalence
\cite[Theorem~3.10]{Orlov2009}
\begin{equation}
\label{eq:factorization-singularity-equivalence}
\MF^\Z(W)\simeq\mathrm D_{\sg}^{\gr}(A).
\end{equation}

Let $X_W=\{W=0\}\subset\PP(w_0,\ldots,w_m)$, and suppose that $X_W$ is smooth. Set
\[
 I_W=\sum_{i=0}^m w_i-d.
\]
If $I_W>0$, \cite[Theorem~2.5(i) and Corollary~2.18]{Orlov2009}
and \eqref{eq:factorization-singularity-equivalence} give an equivalence
\begin{equation}
\label{eq:graded-hypersurface-equivalence}
 \Phi_W\colon\MF^\Z(W)\xrightarrow{\sim}\Ku(X_W),
\end{equation}
where
\[
 \Ku(X_W)= \langle\OO_{X_W},\ldots,\OO_{X_W}(I_W-1)\rangle^\perp.
\]
If $I_W=0$, \cite[Theorem~2.5(iii)]{Orlov2009} gives instead
\begin{equation}
\label{eq:calabi-yau-hypersurface-equivalence}
\MF^\Z(W)\simeq\Db(X_W).
\end{equation}

For an ordinary Fano hypersurface, the equivalence \eqref{eq:graded-hypersurface-equivalence} can be chosen to satisfy the grading compatibility
described in \cite{Toda2013}
\begin{equation}
\label{eq:grading-rotation-compatibility}
  \Phi_W\circ(1)\simeq R_X\circ\Phi_W.
\end{equation}
Then $R_X^d\simeq[2]$. For a weighted Fano hypersurface, we set
\[
 R_{X_W}\coloneq\Phi_W\circ(1)\circ\Phi_W^{-1}.
\]
It is an autoequivalence of $\Ku(X_W)$ satisfying $R_{X_W}^d\simeq[2]$.

\section{Coordinate restrictions on Kuznetsov components}
\label{sec:coordinate-induction}

In this section, we prove the restriction statement used
in hypersurface constructions. Using the graded models in
\eqref{eq:factorization-singularity-equivalence} and
\eqref{eq:graded-hypersurface-equivalence}, we
establish the required phase cuts, mass-Hom bounds and
full numerical support, following
\cite{LLL2026,HLR2025,Cheng2026}.

Let $P=P'[q]$ be a positively graded polynomial ring, where $q$ has weight $a$, and let $W\in P$ have
degree $d$. Write $X$ for its weighted projective hypersurface stack and $Y$ for the
hypersurface defined by $W|_{q=0}$ in the weighted projective stack associated to $P'$. Suppose that
both $X$ and $Y$ are smooth proper Fano stacks. If $I_X$ is the sum of the weights of $P$ minus $d$, their Fano
indices satisfy $I_Y=I_X-a>0$. Consider the following residual categories
\[
 \Ku(X)=\langle\OO_X,\ldots,\OO_X(I_X-1)\rangle^\perp, \qquad
  \Ku(Y)=\langle\OO_Y,\ldots,\OO_Y(I_Y-1)\rangle^\perp.
\]
Applying \eqref{eq:graded-hypersurface-equivalence} to $W$ and $W|_{q=0}$ identifies these
residual categories with the corresponding idempotent-complete graded factorization categories. Let $R_X$
and $R_Y$ denote the corresponding grading shifts.

Set $R=P/(W)$ and $S=R/(q)$. If $q$ is regular in $R$, the adjoint functors
\[
-\otimes_R^{\mathbf L}S \colon\Db(\operatorname{grmod}R)
 \longrightarrow\Db(\operatorname{grmod}S), \qquad \operatorname{Res}_{S/R}
\colon\Db(\operatorname{grmod}S)
\longrightarrow\Db(\operatorname{grmod}R)
\]
preserve perfect complexes and hence induce an adjoint pair
\[
\overline U\colon\mathrm D_{\sg}^{\gr}(R) \longrightarrow\mathrm
  D_{\sg}^{\gr}(S), \qquad \overline V\colon\mathrm D_{\sg}^{\gr}(S) \longrightarrow\mathrm
 D_{\sg}^{\gr}(R).
\]
Via \eqref{eq:graded-hypersurface-equivalence}, we obtain
\[
  U=\Phi_Y\circ\overline U\circ\Phi_X^{-1}, \qquad V=\Phi_X\circ\overline V\circ\Phi_Y^{-1},
\]
so that
\[
 U\colon\Ku(X)\longrightarrow\Ku(Y), \qquad V\colon\Ku(Y)\longrightarrow\Ku(X), \qquad U\dashv V.
\]
The main result of this section is stated as follows.
\begin{proposition}
\label{prop:weighted-hypersurface-restriction}
In the above setting, suppose that $q$ is regular in $R$ and $0<a\leq d/2$. Let
$\widetilde\sigma_X=(\mathcal P_X,Z_X)$ be a numerical stability condition on $\Ku(X)$ with
mass-Hom bounds and
\[
  R_X\widetilde\sigma_X=\widetilde\sigma_X\cdot2/d.
\]
Then
\begin{equation}
 Z_Y=Z_X\circ V_*,\qquad \mathcal P_Y(\phi) =\{E\in\Ku(Y)\mid V(E)\in\mathcal P_X(\phi)\}
 \label{eq-inverse-slicing}
\end{equation}
define a numerical stability condition $\widetilde\sigma_Y$
on $\Ku(Y)$. It has mass-Hom bounds and satisfies the support
property on the full numerical group. Moreover,
\[
R_Y\widetilde\sigma_Y=\widetilde\sigma_Y\cdot2/d.
\]
\end{proposition}

For the rest of this section, we assume the hypotheses of Proposition~\ref{prop:weighted-hypersurface-restriction}. Set $\mathcal B=\Ku(X)$ and
$\mathcal C=\Ku(Y)$, and write $\sigma_{\mathcal B}=\widetilde\sigma_X$, with $\mathcal P_{\mathcal B}=\mathcal P_X$ and $Z_{\mathcal
B}=Z_X$. For the candidate slicing and charge, write $\mathcal P_{\mathcal C}=\mathcal P_Y$ and $Z_{\mathcal C}=Z_Y$. Both
residual categories $\Ku(X)$ and $\Ku(Y)$ are smooth and proper, have classical generators, and have finite-rank
numerical Grothendieck groups, i.e. ranks of $\Knum(\mathcal B)$ and $\Knum(\mathcal C)$ are finite. Fix a
classical generator $B_0$ of $\mathcal B$ and set $G=U(B_0)$.

We use the Ind-completions of the stable enhancements of
$\mathcal B$ and $\mathcal C$. Recall that the Ind-completion $\Ind(\mathcal C)$ consists of filtered colimits of
objects of $\mathcal C$. An object $E\in\Ind(\mathcal C)$
is \emph{compact} if, for every family $\{F_i\}_{i\in I}$,
the natural map
\[
\bigoplus_{i\in I}\Hom(E,F_i) \longrightarrow \Hom\left(E,\bigoplus_{i\in
  I}F_i\right)
\]
is an isomorphism for every family $\{F_i\}_{i\in I}$. The natural embedding
$\mathcal C\hookrightarrow\Ind(\mathcal C)$ identifies
$\mathcal C$ with the full subcategory of compact objects.

The functor $V\colon\mathcal C\to\mathcal B$ extends uniquely to an exact functor preserving
filtered colimits
\[
 \Ind(V)\colon\Ind(\mathcal C)\longrightarrow\Ind(\mathcal B).
\]
We call this its \emph{Ind-extension}. Recall also that $V$ is
\emph{conservative} if it satisfies
\[
 V(E)\simeq0\quad\Longrightarrow\quad E\simeq0.
\]

\begin{lemma}
\label{lem:cut-generation-compactness}
The object $G$ is a classical generator of $\mathcal C$. The Ind-extension
\[
  V\colon\Ind(\mathcal C)\longrightarrow\Ind(\mathcal B)
\]
is conservative, and an object $E\in\Ind(\mathcal C)$ is compact if and only if $V(E)$ is compact.
\end{lemma}

\begin{proof}
We first show that $\mathcal C=\operatorname{thick}(G)$. The
resolution
\[
  0\longrightarrow R(-a)\xrightarrow{q}R \longrightarrow S\longrightarrow0
\]
gives $\operatorname{pd}_R(S)\leq1$. Hence, derived functors $-\otimes_R^{\mathbf L}S$ and $\operatorname{Res}_{S/R}$ are well defined and preserves perfect complexes. Therefore, we have the induced adjoint functors $U$ and $V$
on $\mathrm D_{\sg}^{\gr}(R)$ and $\mathrm D_{\sg}^{\gr}(S)$. Let $F$ be a graded maximal Cohen--Macaulay $S$-module.
Then we have
\[
H^i\bigl(S\otimes_R^{\mathbf L} \operatorname{Res}_{S/R}(F)\bigr) \simeq
 \begin{cases}
F,&i=0,\\
  F(-a),&i=-1,\\
 0,&i\neq0,-1.
 \end{cases}
\]
The truncation triangle induces
\[
F(-a)[1]\longrightarrow UV(F)\longrightarrow F \xrightarrow{\delta}F(-a)[2].
\]
Write
\(
\delta^N\colon F\longrightarrow F(-Na)[2N]
\)
for the composite of $\delta$ and its successive shifts. Taking $N=dj$ and using $(d)\simeq[2]$, we obtain
\(
F(-Na)[2N]\simeq F[2j(d-a)].
\)
Properness of $\mathcal C$ gives
\[
 \sum_{r\in\Z}\dim_\C\Hom_{\mathcal C}(F,F[r])<\infty.
\]
As $d-a>0$, by properness of $\mathcal C$, it follows that
\[
\Hom_{\mathcal C}\bigl(F,F[2j(d-a)]\bigr)=0 \qquad\text{for }j\gg0.
\]
Thus $\delta^N=0$ for some $N$. The octahedral axiom gives
\[
  \operatorname{Cone}(\delta^N) \in \operatorname{thick}\bigl( \operatorname{Cone}(\delta)(-ja)[2j] \mid 0\leq j<N \bigr).
\]
Moreover,
\[
 \begin{aligned}
\operatorname{Cone}(\delta)(-ja)[2j] &\simeq UV(F)(-ja)[2j+1]\\
  &\simeq U\bigl(V(F)(-ja)[2j+1]\bigr) \in\operatorname{thick}(U(B_0)).
 \end{aligned}
\]
Since $\delta^N=0$,
\[
 \operatorname{Cone}(\delta^N)
 \simeq F(-Na)[2N]\oplus F[1].
\]
Consequently, $F\in\operatorname{thick}(G)$. Every object of the idempotent-completed graded
singularity category is a direct summand of a shift of a maximal Cohen--Macaulay module.
Thick closure therefore gives $\mathcal C=\operatorname{thick}(G)$.

We continue to denote the Ind-extensions by $U$ and $V$. Both preserve compact
objects and coproducts, so their adjunction extends to the Ind-completions. For every
$E\in\Ind(\mathcal C)$, we have
\begin{equation}
 \operatorname{RHom}_{\mathcal C}(G,E) \simeq \operatorname{RHom}_{\mathcal
B}(B_0,V(E)).
 \label{eq-generator-adjunction}
\end{equation}
Since $G$ is a compact generator of $\Ind(\mathcal C)$,
\[
V(E)\simeq0 \quad\Longrightarrow\quad \operatorname{RHom}_{\mathcal C}(G,E)\simeq0,
\]
which further gives $ E\simeq0$. Thus $V$ is conservative. Suppose now that $V(E)$ is compact. Then $V(E)\in\mathcal B$, and
properness of $\mathcal B$ implies
\[
\operatorname{RHom}_{\mathcal C}(G,E) \simeq\operatorname{RHom}_{\mathcal B}(B_0,V(E)) \in\Perf(\C).
\]
Set
\(
A=\operatorname{RHom}_{\mathcal C}(G,G)\), and \(N=\operatorname{RHom}_{\mathcal C}(G,E).
\)
Under the dg Morita equivalence
\[
  \Ind(\mathcal C)\simeq\mathrm D(A), \qquad \mathcal C\simeq\Perf(A),
\]
the object $E$ corresponds to the right $A$-module $N$. Due to
smoothness of $A$, we obtain
\[
N\simeq N\otimes_A^{\mathbf L}A \in\operatorname{thick}_{\mathrm D(A)}
 (N\otimes_\C A).
\]
Hence $N\in\Perf(A)$, so $E$ is compact. Conversely, if $E$ is compact, then $E\in\mathcal
C$ and $V(E)\in\mathcal B$, which proves that $V(E)$ is compact.
\end{proof}

\begin{lemma}
\label{lem:cut-slicing}
The slicing and central charge in \eqref{eq-inverse-slicing} define a pre-stability condition on $\mathcal C$.
\end{lemma}

\begin{proof}
Let $E\in\mathcal P_{\mathcal B}(\phi)$. The triangle
\[
  E(-a)\xrightarrow{q}E\longrightarrow VU(E) \longrightarrow E(-a)[1]
\]
and the Gepner relation give
\[
 E(-a)[1]\in \mathcal P_{\mathcal B}\left(\phi+1-\frac{2a}{d}\right).
\]
Since $a\leq d/2$, we obtain
\begin{equation}
 VU\bigl(\mathcal P_{\mathcal B}(\phi)\bigr)
 \subset\mathcal P_{\mathcal B}(\geq\phi).
 \label{eq-phase-monotonicity}
\end{equation}
Fix $\phi\in\R$. Since the compact objects $\mathcal P_{\mathcal B}(\geq\phi)$ are closed under $[1]$,
\cite[Theorem~A.13]{LLL2026} gives a $t$-structure on $\Ind(\mathcal B)$ with aisle
\(
\widehat{\mathcal P}_{\mathcal B}(\geq\phi) =\operatorname{Coprod} \bigl(\mathcal P_{\mathcal
B}(\geq\phi)\bigr).
\)
Here $\operatorname{Coprod}$ denotes closure under coproducts, extensions and direct summands. Write
$\widehat{\mathcal P}_{\mathcal B}(<\phi)$ for its right orthogonal. For $B\in\mathcal B$, the original
phase truncation triangle
\[
B_{\geq\phi}\longrightarrow B\longrightarrow B_{<\phi}
\]
is also a truncation triangle in $\Ind(\mathcal B)$. Its first term lies in the generated aisle, and phase orthogonality places its last term in
the right orthogonal, since this orthogonality is preserved under coproducts, extensions and direct summands in the first
argument. Thus the extended phase cuts restrict to the original ones on compact objects. On $\Ind(\mathcal C)$, consider the compactly
generated t-structure with aisle
\[
 \widehat{\mathcal D}_{\mathcal C,\geq\phi} = \operatorname{Coprod}
\bigl(U(\mathcal P_{\mathcal B}(\geq\phi))\bigr).
\]
By adjunction, its orthogonal part is given by
\begin{equation}
\widehat{\mathcal D}_{\mathcal C,<\phi}= \bigl\{ E\mid V(E)\in \widehat{\mathcal P}_{\mathcal B}(<\phi) \bigr\}.
 \label{eq-ind-coaisle}
\end{equation}
Since $V$ is exact and preserves coproducts,
\eqref{eq-phase-monotonicity} also gives
\(
V\bigl(\widehat{\mathcal D}_{\mathcal C,\geq\phi}\bigr)
\subset\widehat{\mathcal P}_{\mathcal B}(\geq\phi)\). Thus $V$ is t-exact for
these phase cuts.

Let $E\in\mathcal C$. Applying $V$ to the truncation triangle of $E$ gives the following triangle
\[
(VE)_{\geq\phi}\longrightarrow VE\longrightarrow(VE)_{<\phi}.
\]
Since $VE\in\mathcal B$, both $(VE)_{\geq\phi}$ and $(VE)_{<\phi}$ are compact. By
Lemma~\ref{lem:cut-generation-compactness}, $E_{\geq\phi}$ and $E_{<\phi}$ are
compact as well. Hence the phase-cut t-structure restricts to $\mathcal C$.

Moreover, conservativity gives
\[
 \begin{aligned}
 VE\in\mathcal P_{\mathcal B}(\geq\phi) &\Longrightarrow V(E_{<\phi})=0 \Longrightarrow E_{<\phi}=0,\\
VE\in\mathcal P_{\mathcal B}(<\phi) &\Longrightarrow V(E_{\geq\phi})=0
  \Longrightarrow E_{\geq\phi}=0.
 \end{aligned}
\]
Consequently,
\[
 \mathcal D_{\mathcal C,\geq\phi} =V^{-1}\bigl(\mathcal P_{\mathcal B}(\geq\phi)\bigr),\qquad \mathcal D_{\mathcal C,<\phi}
  =V^{-1}\bigl(\mathcal P_{\mathcal B}(<\phi)\bigr).
\]
Their boundedness follows from that of the base phase cuts. We now verify that
\[
\mathcal P_{\mathcal C}(\phi) = \{E\in\mathcal C\mid V(E)\in\mathcal P_{\mathcal B}(\phi)\}
\]
defines a slicing. Exactness of $V$ gives
\(
\mathcal P_{\mathcal C}(\phi+1) = \mathcal P_{\mathcal C}(\phi)[1]\). If $\phi_1>\phi_2$,
then
\[
\mathcal P_{\mathcal C}(\phi_1) \subset\mathcal D_{\mathcal C,\geq\phi_1}, \qquad \mathcal P_{\mathcal C}(\phi_2) \subset\mathcal
D_{\mathcal C,<\phi_1}.
\]
Orthogonality therefore implies
\(
\Hom_{\mathcal C} \bigl(\mathcal P_{\mathcal C}(\phi_1), \mathcal P_{\mathcal C}(\phi_2)\bigr)=0.
\)
For the HN property, let $B^{(0)}=V(E)$ and write its successive highest-phase truncation triangles as
\[
 F_i\longrightarrow B^{(i-1)} \longrightarrow B^{(i)}, \qquad F_i\in\mathcal P_{\mathcal
B}(\phi_i), \qquad \phi_1>\cdots>\phi_r, \qquad B^{(r)}=0.
\]
Starting with $E^{(0)}=E$, truncate successively at $\phi_i$. T-exactness gives triangles
\[
A_i\longrightarrow E^{(i-1)} \longrightarrow E^{(i)}
\]
with
\(
V(A_i)\simeq F_i,\) and \(V(E^{(i)})\simeq B^{(i)}. \) Hence $A_i\in\mathcal P_{\mathcal C}(\phi_i)$, and conservativity gives
\[
V(E^{(r)})=0\quad\Longrightarrow\quad E^{(r)}=0.
\]
These triangles form a finite HN filtration of $E$, whose image is
the filtration of $V(E)$.

For $0\neq E\in\mathcal P_{\mathcal C}(\phi)$, conservativity ensures $V(E)\neq0$. Thus
\[
Z_{\mathcal C}(E) = Z_{\mathcal B}(V(E)) \in\R_{>0}e^{i\pi\phi}.
\]
The identification of the HN factors also gives
\begin{equation}
 \begin{aligned}
  m_{\widetilde\sigma_Y}(E) &=\sum_{i=1}^r|Z_{\mathcal B}(F_i)|
 =m_{\widetilde\sigma_X}(V(E)).
 \end{aligned}
 \label{eq-mass-identity}
\end{equation}
\end{proof}

\begin{lemma}
\label{lem:cut-mass-support}
The pre-stability condition of Lemma~\ref{lem:cut-slicing} is numerical and satisfies the support property on $\Knum(\mathcal
C)$. The resulting stability condition has mass-Hom bounds.
\end{lemma}

\begin{proof}
 If $v(E)=0$ in $\Knum(\mathcal C)$, adjunction gives
\[
  \chi_{\mathcal B}(B,V(E)) =\chi_{\mathcal C}(U(B),E)=0 \qquad\text{for all }B\in\mathcal
 B.
\]
Thus $v(V(E))=0$ in $\Knum(\mathcal B)$. The map $V_*$ then descends to the numerical groups, and
$Z_{\mathcal C}$ is numerical. Choose $A_1,\ldots,A_s$ such that $\{\chi_{\mathcal C}(A_j,-)\}_{j=1}^s$
spans $\Knum(\mathcal C)^{\vee}_{\mathbb R}$. For any fixed norm there is a constant $c$ such that
\begin{equation}
  \|v(E)\|\leq c\sum_{j=1}^s|\chi_{\mathcal C}(A_j,E)|.
 \label{eq-euler-norm}
\end{equation}
Every $A_j$ belongs to the thick subcategory generated by $G=U(B_0)$. A fixed finite construction of $A_j$ from $G$ by cones, shifts, sums and
direct summands yields a constant $c_j$ with
\[
 \dim\operatorname{Hom}_{\mathcal C}^*(A_j,E)
 \leq c_j\dim\operatorname{Hom}_{\mathcal C}^*(G,E).
\]
By adjunction and Lemma~\ref{lem-total-mass-hom}, the last
quantity is at most
\(
C_{B_0}m_{\sigma_{\mathcal B}}(V(E)). \) Combining this with \eqref{eq-euler-norm} proves
\begin{equation}
\|v(E)\|\leq C m_{\sigma_{\mathcal B}}(V(E)) =C m_{\widetilde\sigma_Y}(E).
 \label{eq-full-numerical-bound}
\end{equation}
For every $\widetilde\sigma_Y$-semistable object $E$, this becomes
\(
\|v_{\mathcal C}(E)\| \leq C|Z_{\mathcal C}(E)|, \) which is the support property on the
full numerical Grothendieck group. Applying the same thick-generation estimate to
an arbitrary fixed object $A$ instead of the finite collection
$A_j$ proves a mass-Hom bound on $\mathcal C$, as in
\cite[Lemma~2.8]{HLR2025} and \cite[Lemma~7.4]{LLL2026}.
\end{proof}

\begin{proof}[Proof of Proposition~\ref{prop:weighted-hypersurface-restriction}]
Lemma~\ref{lem:cut-generation-compactness} supplies the generator and compactness properties
required in Lemma~\ref{lem:cut-slicing}. The latter gives the pre-stability condition
$\widetilde\sigma_Y$ in \eqref{eq-inverse-slicing}.
Lemma~\ref{lem:cut-mass-support} proves its
numericality, mass-Hom bounds, and full numerical support. Since $V$
commutes with the grading shifts, we have
\(
V R_Y\simeq R_X V. \) For every $\phi\in\mathbb R$, the source Gepner equation and
\eqref{eq-inverse-slicing} therefore give
\[
R_Y\bigl(\mathcal P_Y(\phi)\bigr) =\mathcal P_Y\left(\phi+\frac2d\right).
\]
Moreover,
\[
Z_Y(R_Y E)=Z_X(R_XV(E)) =e^{2\pi\ii/d}Z_Y(E).
\]
This concludes that $R_Y\widetilde\sigma_Y=\widetilde\sigma_Y\cdot2/d$.
\end{proof}

\section{Gepner-type stability from elliptic curves}
\label{sec:elliptic-seeds}

In this section, we construct the  stability conditions on kuznetsov components of Fano hypersurfaces of degrees three, four and six from products of elliptic curves, based on
\cite{LMP2025,LLL2026,BFK2014}.

For $d\in\{3,4,6\}$, let $w_d$ be the polynomial in the following table. The line bundle $\mathcal H_d=\mathcal O_{E_d}(1)$
on the weighted projective curve $E_d=\{w_d=0\}$ has degree $c_d$.
\[
\begin{array}{c|c|c|c|c}
 d&w_d&\text{weights}&c_d&\tau_d\\ \hline
 3&x^3+y^3+z^3&(1,1,1)&3&-\frac32+i\frac{\sqrt3}{2}\\
 4&x^4+y^4+z^2&(1,1,2)&2&-1+i\\
 6&x^6+y^3+z^2&(1,2,3)&1&-\frac12+i\frac{\sqrt3}{2}
\end{array}
\]
These curves are smooth elliptic curves without stacky points. For the
quartic curve, the only point of $\Pp(1,1,2)$ with nontrivial stabilizer
does not lie on the curve. On the sextic curve, a point with $x=0$ has
both $y$ and $z$ nonzero. Their relatively prime weights give a trivial
stabilizer. Every other point has a nonzero coordinate of weight one.
Adjunction gives a trivial canonical bundle in each case. The degree of
the polarization is $d$ divided by the product of the three weights,
which gives the values of $c_d$.

Fix $d$, abbreviate $E=E_d$, $c=c_d$ and $\tau=\tau_d$, and write
\[
\lambda=\exp(2\pi i/d)=1+\tau, \qquad \tau^2+c\tau+c=0.
\]
For an integer $k\geq2$, let
\[
W_{d,k}=\sum_{j=1}^k w_d(x_j,y_j,z_j).
\]
Let $X_{d,k}$ be its weighted projective hypersurface, which is a smooth proper Fano stack, since the sum
of the weights is $kd>d$. By \eqref{eq:graded-hypersurface-equivalence}, $\Ku(X_{d,k})$ is
equivalent to $\MF^\Z(W_{d,k})$. We use the idempotent-complete enhancement of this graded model
in the proof. The main result of this section is stated as follows.

\begin{theorem}\label{thm:elliptic-block-seed}
For every $d\in\{3,4,6\}$ and $k\geq2$, 
$\Ku(X_{d,k})$ admits a numerical stability condition
$\widetilde\sigma$ with support on its full numerical Grothendieck
group. It satisfies
\[
  R_{X_{d,k}}\widetilde\sigma=\widetilde\sigma\cdot\frac2d.
\]
For every fixed object $P$, there is a constant $C_P>0$ such that
\[
\sum_{a\in\Z}\dim_{\C}\Hom(P,F[a]) \leq C_Pm_{\widetilde\sigma}(F) \qquad\text{for all }F\in\Ku(X_{d,k}).
\]
\end{theorem}

Let
\[
L=\frac{\bigoplus_{j=1}^k\Z e_j} {\langle d(e_j-e_k)\mid 1\leq j<k\rangle}.
\]
Give the variables in the $j$-th block their weights times $e_j$. Every block potential has degree $de_j$, and these degrees are equal in
$L$. The degree map $L\to\Z$ sending every $e_j$ to one has kernel
\[
\Gamma=\bigoplus_{j=1}^{k-1}(\Z/d)(e_j-e_k).
\]

\begin{lemma}
\label{lem:elliptic-grading-comparison}
There are equivalences of enhanced categories
\[
 \MF^L(W_{d,k})\simeq\Db(E^k),
 \qquad
 \MF^\Z(W_{d,k})\simeq\Db(E^k)^\Gamma.
\]
Under the first equivalence, the shift $e_j$ acts by the transform
\[
  T=\operatorname{ST}_{\OO_E}\circ (-\otimes\mathcal H_d)
\]
on the $j$-th factor. Denote this transform by $T_j$. The action of
$e_j-e_k\in\Gamma$ is $T_jT_k^{-1}$. Under the second equivalence, the
ordinary grading shift $(1)$ has underlying transform
$T_k$.
\end{lemma}

\begin{proof}
By \cite[Corollary~2.40 and Theorem~3.5]{BFK2014} and
\eqref{eq:calabi-yau-hypersurface-equivalence}, there
is an equivalence of enhanced categories
\begin{equation}
  \Phi\colon\MF^L(W_{d,k})\xrightarrow{\sim}\Db(E^k).
 \label{eq:elliptic-fine-category}
\end{equation}
For each factor, the grading shift corresponds to
\[
T=\operatorname{ST}_{\mathcal O_E} \circ(-\otimes\mathcal H_d), \qquad T^d\simeq[2],
\]
by \cite[Proposition~3.2]{Toda2013}. Compatibility with the factor equivalences gives
\(
\Phi\circ(e_j)\simeq T_j\circ\Phi.
\)
Consider the diagonalizable group with character group $L$
\[
G_L= \{(t_1,\ldots,t_k)\in(\C^*)^k \mid t_1^d=\cdots=t_k^d\}.
\]
Its diagonal subgroup is $\C^*$ and its quotient is $Q=(\mu_d)^{k-1}$. Consequently the
finely graded category is the $Q$-equivariant version of the ordinarily graded
category. Since $Q^\vee\simeq\Gamma$, finite abelian equivariantization duality
\cite[Theorem~4.2]{Elagin2015} gives
\begin{equation}\label{eq:elliptic-coarse-category}
 \MF^{\Z}(W_{d,k})\simeq\Db(E^k)^\Gamma.
\end{equation}

The character action of $\gamma\in\Gamma$ on $\MF^L(W_{d,k})$ is the
grading shift $(\gamma)$. In particular,
\[
\Phi\circ(e_j-e_k) \simeq T_jT_k^{-1}\circ\Phi.
\]
The composition isomorphisms of grading shifts give
\[
  (\gamma)\circ(\gamma')\simeq(\gamma+\gamma'), \qquad (0)\simeq\operatorname{id},
\]
and hence the coherence isomorphisms for the transported
$\Gamma$-action. The relation $d(e_j-e_k)=0$ gives
\(
(T_jT_k^{-1})^d\simeq\operatorname{id}
\)
with this coherence. Since $e_k$ maps to $1\in\Z$, the ordinary grading shift
corresponds to the autoequivalence induced by $(e_k)$. If
\(
 \operatorname{For}\colon\Db(E^k)^\Gamma\longrightarrow\Db(E^k)
\)
denotes the forgetful functor, then, under
\eqref{eq:elliptic-coarse-category},
\(
\operatorname{For}\circ(1) \simeq T_k\circ\operatorname{For}.
\)
This proves the assertions.
\end{proof}

Choose an origin on $E$, let $h_j$ be its degree-one class on the $j$-th factor of $A=E^k$, and put $h=\sum_jh_j$. We then set
\begin{equation}\label{eq:elliptic-product-charge}
 Z_A(F)=-\int_A e^{-\tau h}\ch(F).
\end{equation}

Let $\Lambda_k=\Z^{\{I\subseteq\{1,\ldots,k\}\}}$ and define
\[
  v_k(F)=\left(\int_A h^I\ch_{k-|I|}(F)\right)_I, \qquad h^I=\prod_{i\in
   I}h_i, \quad h^\varnothing=1.
\]
These coordinates are integral by \cite[Lemma~4.3]{LMP2025}.

\begin{lemma}
\label{lem:elliptic-product-gepner}
There is a stability condition $\widetilde\sigma_A$ on $\Db(A)$
with respect to $(\Lambda_k,v_k)$, with central charge
\eqref{eq:elliptic-product-charge} and with every point sheaf stable
of phase one. It satisfies
\[
  T_j\widetilde\sigma_A =\widetilde\sigma_A\cdot\frac2d \text{ for }1\leq j\leq k.
\]
In particular, the $\Gamma$-action preserves $\widetilde\sigma_A$.
\end{lemma}

\begin{proof}
By \cite[Proposition~4.5]{LLL2026}, there is a stability condition
$\widetilde\sigma_A=(\mathcal P_A,Z_A)$ with respect to
$(\Lambda_k,v_k)$, with central charge
\eqref{eq:elliptic-product-charge}, for which every point
sheaf is stable of phase one. By
\cite[Theorem~3.11]{LMP2025}, its charge and the
phases of the point sheaves determine its slicing.

Fix $j$ and let $p\colon A\to E^{k-1}$ omit the $j$-th factor. For
$F\in\Db(A)$, write
\[
 a_F=p_*\bigl(h_j\ch(F)\bigr), \qquad b_F=p_*\ch(F).
\]
Tensoring by the polarization on this factor sends
\(
(a_F,b_F)\) to \((a_F,b_F+ca_F).\) The evaluation triangle for the partial spherical twist gives
\[
  [T_jF] = [F\otimes\mathcal H_{d,j}] - [p^*\mathbf Rp_*(F\otimes\mathcal H_{d,j})].
\]
Since $\operatorname{td}(E)=1$, Grothendieck--Riemann--Roch therefore yields
\[
 \begin{pmatrix}a_{T_jF}\\ b_{T_jF}\end{pmatrix}
 =
 \begin{pmatrix}1-c&-1\\ c&1\end{pmatrix}
 \begin{pmatrix}a_F\\ b_F\end{pmatrix}.
\]
Using $\tau^2+c\tau+c=0$ and $\lambda=1+\tau$, we obtain
\[
 \begin{aligned}
b_{T_jF}-\tau a_{T_jF} &= \bigl(c-\tau(1-c)\bigr)a_F+(1+\tau)b_F\\
 &=\lambda(b_F-\tau a_F).
 \end{aligned}
\]
Set $h_{\widehat j}=\sum_{i\ne j}h_i$. Then
\[
Z_A(F) = -\int_{E^{k-1}} e^{-\tau h_{\widehat j}}(b_F-\tau a_F),
\]
and hence
\begin{equation}
Z_A(T_jF)=\lambda Z_A(F).
 \label{eq:elliptic-eigencharge}
\end{equation}
Components of degree one on the $j$-th factor contribute neither to $a_F$ nor to $b_F$, so
the calculation also includes mixed cohomology classes.

We next determine the phase lift. On $E$, the numerical action of $T$ is
\[
T_*(r,e)=\bigl((1-c)r-e,cr+e\bigr).
\]
For $x=(x_1,\ldots,x_k)$, set
\[
  H_x=T_j^{-1}\mathcal O_x = L_{x_j}\boxtimes \mathop{\boxtimes}_{i\ne j}\mathcal O_{x_i}, \qquad \deg L_{x_j}=1-c.
\]
Partial Fourier--Mukai transforms on the factors $i\ne j$ take $H_x$ to an
exterior product of line bundles. This is a simple semihomogeneous vector bundle and is
stable for every numerical stability condition \cite[Corollary~2.16]{FLZ2022}. Applying this to the
transported condition shows that $H_x$ is $\widetilde\sigma_A$-stable.

Moreover,
\[
\Hom(H_x,\mathcal O_x)\ne0, \qquad \Hom(\mathcal O_x,H_x[1])\ne0.
\]
Since $\phi_A(\mathcal O_x)=1$, the phase vanishing axiom implies
\(
0\leq\phi_A(H_x)\leq1. \) On the other hand,
\[
Z_A(H_x) =\tau+c-1 =-\lambda^{-1} =\exp\left(i\pi\left(1-\frac2d\right)\right).
\]
Thus
\(
\phi_A(H_x)=1-2/d. \) Consequently, every point sheaf is stable of phase $1-2/d$ for both $T_j\widetilde\sigma_A$ and
$\widetilde\sigma_A\cdot2/d$. Their central charges also agree, since
\[
Z_{T_j\widetilde\sigma_A} =Z_A\circ T_j^{-1} =\lambda^{-1}Z_A =Z_{\widetilde\sigma_A\cdot2/d}.
\]
By \cite[Theorem~3.11]{LMP2025},
\begin{equation}
T_j\widetilde\sigma_A =\widetilde\sigma_A\cdot\frac2d \text{ for }1\leq j\leq
k.
 \label{eq:elliptic-factor-gepner}
\end{equation}
In particular, for $j<k$,
\(
(T_jT_k^{-1})\widetilde\sigma_A =\widetilde\sigma_A,
\)
so the coherent $\Gamma$-action preserves $\widetilde\sigma_A$.
\end{proof}

We fix the stability condition $\widetilde\sigma_A=(\mathcal P_A,Z_A)$ of
Lemma~\ref{lem:elliptic-product-gepner} for the remainder of
this section.

\begin{lemma}
\label{lem:elliptic-twist-mass}
For every subset $I\subset\{1,\ldots,k\}$, there is a constant $C_I>0$ such that
\[
m_{\widetilde\sigma_A} \bigl(F\otimes\OO_A(-4h_I)\bigr) \leq C_I m_{\widetilde\sigma_A}(F) \qquad\text{ for all }F\in\Db(A), \qquad
 h_I=\sum_{j\in I}h_j.
\]
\end{lemma}

\begin{proof}
Write $\tau=B+i\omega$. In each case,
\(
B,\omega^2\in\Q,\) and \( Z_A\bigl(K_0(A)\bigr)\subset\Q+i\omega\Q.
\)
Fix $j$ and perform the product construction with the $j$-th factor last. On the remaining factors, apply
\[
S\colon\C\longrightarrow\C, \qquad S(x+iy)=x+i\frac{y}{\omega}.
\]
The transformed central charge takes values in $\Q+i\Q$, so its heart is noetherian. If the original input
polynomial is
\(
am+b+i(cm+e),
\)
the transformed polynomial is
\[
am+b+i\left(\frac c\omega m+\frac e\omega\right).
\]

For this fixed input, \cite[Theorem~4.1]{LMP2025}, with
$s=\omega^2$ and $t=1$, gives a stability condition for
every $\beta\in\R$. Applying $S^{-1}$ to its charge gives
\[
 \begin{aligned}
b-\beta a+\omega c+i(e-\omega a-\beta c) &=(b+ie)-(\beta+i\omega)(a+ic).
 \end{aligned}
\]
Thus the resulting condition $\widetilde\sigma_\beta$ has charge
\[
 Z_\beta(F)
 =
 -\int_A
 e^{-\tau\sum_{i\ne j}h_i-(\beta+i\omega)h_j}\ch(F).
\]
All these conditions use the same numerical quotient lattice
$\Lambda$, determined by the fixed input. Moreover,
\(
Z_\beta(\mathcal O_x)=-1\), and
 \(
\mathcal O_x\in\mathcal P_\beta(1)\) for every $x\in A$. By \cite[Theorem~3.11]{LMP2025}, we obtain
\(
\widetilde\sigma_B=\widetilde\sigma_A.
\)

We claim that
\[
[B,B+4]\longrightarrow\Stab_\Lambda(A), \qquad \beta\longmapsto\widetilde\sigma_\beta
\]
is continuous. Fix $\beta_0$. By \cite[Theorem~7.1]{Bridgeland}, the charges $Z_\beta$ admit a local
continuous lift $\widetilde\sigma'_\beta$ with
\(
 \widetilde\sigma'_{\beta_0}
 =\widetilde\sigma_{\beta_0}.
\) Every skyscraper sheaf remains stable \cite{FLZ2022}. Since
$Z_\beta(\mathcal O_x)=-1$, continuity of its phase gives
\(
\phi_{\widetilde\sigma'_\beta}(\mathcal O_x)=1.
\)
Hence \cite[Theorem~3.11]{LMP2025} implies
\(
\widetilde\sigma'_\beta=\widetilde\sigma_\beta
\)
near $\beta_0$, proving the claim.

At $[B+4]$, we have
\(
Z_{B+4}(F)=Z_A(F\otimes \mathcal O_A(-4h_j)).
\)
Both $\widetilde\sigma_{B+4}$ and $\widetilde\sigma_A\otimes \mathcal
O_A(4h_j)$ have point sheaves stable of phase one. Uniqueness therefore gives
\(
 \widetilde\sigma_{B+4}
 =\widetilde\sigma_A\otimes \mathcal O_A(4h_j).
\)
The continuous path lies in a single component of finite Bridgeland distance. Consequently,
\[
 D_j = d\bigl(\widetilde\sigma_A, \widetilde\sigma_A\otimes \mathcal O_A(4h_j)\bigr)
<\infty.
\]
Thus,
\[
 \begin{aligned}
 m_{\widetilde\sigma_A}(F\otimes \mathcal O_A(-4h_j))
 &=
 m_{\widetilde\sigma_A\otimes \mathcal O_A(4h_j)}(F)\\
 &\leq e^{D_j}m_{\widetilde\sigma_A}(F).
 \end{aligned}
\]
Applying this estimate successively for $j\in I$ gives
\[
 m_{\widetilde\sigma_A} \bigl(F\otimes\mathcal O_A(-4h_I)\bigr) \leq \exp\left(\sum_{j\in
I}D_j\right) m_{\widetilde\sigma_A}(F).
\]
Thus one may take $C_I=\exp(\sum_{j\in I}D_j)$.
\end{proof}

\begin{lemma}
\label{lem:elliptic-product-mass-hom}
The product condition $\widetilde\sigma_A$ has mass-Hom bounds.
More precisely, for every $P\in\Db(A)$, there is a constant
$C_P>0$ such that
\[
\sum_{a\in\Z}\dim_\C\Hom(P,F[a]) \leq C_Pm_{\widetilde\sigma_A}(F) \qquad\text{for all }F\in\Db(A).
\]
\end{lemma}

\begin{proof}
Let
\[
 f=\prod_{j=1}^k f_j\colon A\longrightarrow Y=(\Pp^1)^k
\]
be the product of the elliptic involution quotients. By \cite[Proposition~5.1]{LLL2026}, the formulas
\[
\mathcal P_Y(\phi) = \{E\in\Db(Y)\mid f^*E\in\mathcal P_A(\phi)\}, \qquad Z_Y=Z_A\circ f^*
\]
define a stability condition $\widetilde\sigma_Y$. Write $\epsilon_j$ for the involution on the $j$-th
factor. Since $\epsilon_j^*h_i=h_i$ for all $i$, charge and point-phase uniqueness gives
\(
\epsilon_j^*\widetilde\sigma_A=\widetilde\sigma_A.
\)
By \cite[Propositions~4.5 and~5.1]{LLL2026}, the two-parameter family
containing $\widetilde\sigma_A$ induces a continuous family on $Y$, containing an
\emph{algebraic stability condition}, i.e. a stability condition
whose heart is of finite length with finitely many isomorphism classes of
simple objects, up to the
$\widetilde{\mathrm{GL}}_2^+(\R)$-action
\cite[Theorem~5.4]{LLL2026}.

Let
$\widetilde\sigma_{\mathrm{alg}}$ be an algebraic stability condition,
with heart $\mathcal H$ and simples $S_1,\ldots,S_N$.
Their phases lie in an interval of width less than one,
so there is a real linear functional $\ell$ satisfying
\[
\eta=\min_j\ell(Z_{\mathrm{alg}}(S_j))>0.
\]
For $E\in\mathcal H$, a composition series gives
\[
\eta\,\operatorname{length}_{\mathcal H}(E) \leq\ell(Z_{\mathrm{alg}}(E))
   \leq\|\ell\|\,|Z_{\mathrm{alg}}(E)|.
\]
For a fixed $P\in\Db(Y)$, set
\(
M_P=\max_j \sum_{a\in\Z}\dim_\C\Hom_{Y}(P,S_j[a])<\infty.
\)
The long exact sequences associated to the composition series yield
\[
\sum_{a\in\Z}\dim_\C\Hom_{Y}(P,E[a]) \leq M_P\operatorname{length}_{\mathcal H}(E) \leq\frac{M_P\|\ell\|}{\eta}|Z_{\mathrm{alg}}(E)|.
\]
Shifting a semistable object into $\mathcal H$ and then applying this estimate to the HN factors gives
\[
 \sum_{a\in\Z}\dim_\C\Hom_{Y}(P,E[a])\leq C_Pm_{\widetilde\sigma_{\mathrm{alg}}}(E),
\]
which is preserved by the $\widetilde{\mathrm{GL}}_2^+(\R)$-action and by
finite-distance mass comparison \cite[Corollary~2.14]{HLR2025}. The continuous
family connects the algebraic stability condition to
$\widetilde\sigma_Y$, so
\[
  \sum_{a\in\Z}\dim_\C\Hom_{Y}(P,E[a])\leq C_Pm_{\widetilde\sigma_Y}(E),
\]
for any $E\in\Db(Y)$, as in \cite[Remark~7.3]{LLL2026}. For the partial quotient on the $j$-th factor, the kernel
filtration gives a triangle
\[
 \epsilon_j^*F\otimes\mathcal O_A(-4h_j)
 \longrightarrow f_j^*f_{j*}F
 \longrightarrow F
 \longrightarrow
 \epsilon_j^*F\otimes\mathcal O_A(-4h_j)[1].
\]
Iterating these triangles gives a finite filtration of $f^*f_*F$
with factors
\begin{equation}
F_I=\epsilon_I^*F\otimes\mathcal O_A(-4h_I), \qquad \epsilon_I=\prod_{j\in
I}\epsilon_j, \quad h_I=\sum_{j\in I}h_j, \quad I\subset\{1,\ldots,k\},
 \label{eq:elliptic-branch-factors}
\end{equation}
by \cite[Example~3.20 and Proposition~3.21]{LLL2026}. Lemma~\ref{lem:elliptic-twist-mass}
and invariance under the involutions give constants $C_I$ such that
\(
m_{\widetilde\sigma_A}(F_I) \leq C_I m_{\widetilde\sigma_A}(F). \) Thus,
\begin{equation}
 \begin{aligned}
 m_{\widetilde\sigma_Y}(f_*F)
 &=m_{\widetilde\sigma_A}(f^*f_*F)\\
  &\leq\sum_I m_{\widetilde\sigma_A}(F_I)\\
 &\leq C m_{\widetilde\sigma_A}(F),
 \end{aligned}
 \label{eq:elliptic-reverse-mass}
\end{equation}
where $C=\sum_I C_I.$ Choose a classical generator $G$ of $\Db(Y)$. For $F\in
D_{\mathrm{qc}}(A)$, adjunction gives
\[
 \RHom_A(f^*G,F)\simeq\RHom_Y(G,f_*F).
\]
Consequently,
\[
\RHom_A(f^*G,F)=0 \Longrightarrow f_*F=0.
\]
Due to conservativity of finite pushforward, we obtain $F=0$. Hence $f^*G$ is a compact generator of $D_{\mathrm{qc}}(A)$,
and
\(
\operatorname{thick}(f^*G) =\operatorname{Perf}(A)=\Db(A).
\)

Finally, adjunction and \eqref{eq:elliptic-reverse-mass} imply
\[
 \begin{aligned}
 \sum_{a\in\Z}\dim_\C\Hom_{A}(f^*G,F[a]) &=\sum_{a\in\Z}\dim_\C\Hom_{Y}(G,f_*F[a])\\
  &\leq C_Gm_{\widetilde\sigma_Y}(f_*F)\\
 &\leq CC_Gm_{\widetilde\sigma_A}(F).
 \end{aligned}
\]
For every fixed $P\in\operatorname{thick}(f^*G)$, the long exact sequence of
morphisms gives a constant $c_P$ such that
\[
\sum_{a\in\Z}\dim_\C\Hom_{A}(P,F[a])\leq
c_P\sum_{a\in\Z}\dim_\C\Hom_{A}(f^*G,F[a]).
\]
Therefore
\[
\sum_{a\in\Z}\dim_\C\Hom_{A}(P,F[a]) \leq c_PCC_Gm_{\widetilde\sigma_A}(F).
\]
This proves the total mass--Hom bound.
\end{proof}

Set $\mathcal C=\Db(A)^\Gamma$, and let
$V\colon\mathcal C\longrightarrow\Db(A)$ be the forgetful functor.

\begin{lemma}
\label{lem:elliptic-equivariant-support}
The charge and slicing
\[
 Z_{\mathcal C}=Z_A V_* ,\qquad \mathcal P_{\mathcal C}(\phi)=V^{-1}\mathcal P_A(\phi)
\]
define a numerical stability condition on $\mathcal C$ with support
on its full numerical Grothendieck group and mass-Hom bounds.
\end{lemma}

\begin{proof}
By \eqref{eq:elliptic-factor-gepner}, the $\Gamma$-action preserves $\widetilde\sigma_A$. Exactness of $V$ gives
\(
\mathcal P_{\mathcal C}(\phi+1) =\mathcal P_{\mathcal C}(\phi)[1]. \) For $F_i\in\mathcal P_{\mathcal C}(\phi_i)$ with $\phi_1>\phi_2$,
\[
\Hom_{\mathcal C}(F_1,F_2) \hookrightarrow\Hom_A(VF_1,VF_2)=0.
\]
Moreover, conservativity gives
\[
 0\ne F\in\mathcal P_{\mathcal C}(\phi)
 \Longrightarrow
 Z_{\mathcal C}(F)
 =Z_A(VF)\in\R_{>0}e^{i\pi\phi}.
\]

Let
\[
 0=E_0\longrightarrow E_1\longrightarrow\cdots
 \longrightarrow E_r=VF
\]
be the HN filtration, with factors $A_i\in\mathcal P_A(\phi_i)$ and
$\phi_1>\cdots>\phi_r$. Its uniqueness equips the filtration with compatible
$\Gamma$-equivariant structures. It therefore lifts to
\[
 0=F_0\longrightarrow F_1\longrightarrow\cdots \longrightarrow F_r=F, \qquad
 V(F_i/F_{i-1})\simeq A_i.
\]
Thus $\mathcal P_{\mathcal C}$ has the HN property, and
\begin{equation}
m_{\mathcal C}(F) =\sum_i|Z_A(A_i)| =m_{\widetilde\sigma_A}(VF).
 \label{eq:elliptic-descent-mass}
\end{equation}

For a fixed $P\in\mathcal C$, the inclusions
\[
 \Hom_{\mathcal C}(P,F[a])
 \hookrightarrow\Hom_A(VP,VF[a])
\]
and Lemma~\ref{lem:elliptic-product-mass-hom} give
\[
 \begin{aligned}
\sum_a\dim\Hom_{\mathcal C}(P,F[a]) &\leq\sum_a\dim\Hom_A(VP,VF[a])\\
&\leq C_{VP}m_{\widetilde\sigma_A}(VF)\\
  &=C_{VP}m_{\mathcal C}(F).
 \end{aligned}
\]
Hence the descended pre-stability condition has total mass--Hom bounds. Let
$I\colon\Db(A)\to\mathcal C$ be equivariant induction, so
that $I\dashv V$. Adjunction gives
\(
\chi_A(B,VF)=\chi_{\mathcal C}(IB,F),
\)
for $B\in\Db(A)$. Consequently, if $v_{\mathcal C}(F)=0$, then $v_A(VF)=0$. Thus $V_*$
descends to a homomorphism
\(
V_{\mathrm{num}}\colon \Knum(\mathcal C)\to\Knum(A),
\)
and $Z_{\mathcal C}$ is numerical.

The category $\mathcal C$ is smooth and proper, and $\Knum(\mathcal C)$ has finite rank by
\eqref{eq:elliptic-coarse-category} and \eqref{eq:graded-hypersurface-equivalence}. Choose $P_1,\ldots,P_N\in\mathcal C$ such
that
\[
v\longmapsto \bigl(\chi_{\mathcal C}(P_1,v),\ldots, \chi_{\mathcal C}(P_N,v)\bigr)
\]
is injective on $\Knum(\mathcal C)_\R$. Define
\[
\|v\|^2 =\sum_{j=1}^N\chi_{\mathcal C}(P_j,v)^2.
\]
For a semistable object $F$, the property of mass--Hom bounds yields
\[
|\chi_{\mathcal C}(P_j,F)| \leq\sum_a\dim\Hom_{\mathcal C}(P_j,F[a]) \leq C_j|Z_{\mathcal C}(F)|.
\]
Taking $C=(\sum_jC_j^2)^{1/2}$, we obtain
\[
\|v_{\mathcal C}(F)\| \leq C|Z_{\mathcal C}(F)|.
\]
Therefore the quadratic form
\[
 Q(v)=C^2|Z_{\mathcal C}(v)|^2-\|v\|^2
\]
satisfies
\[
 Q(v_{\mathcal C}(F))\geq0
 \quad\text{for every semistable }F,
 \qquad
 Q|_{\ker Z_{\mathcal C}}=-\|\cdot\|^2.
\]
This proves the support property on the full numerical group,
following the argument in \cite[Theorem~2.1]{Cheng2026}.
\end{proof}

\begin{proof}[Proof of Theorem~\ref{thm:elliptic-block-seed}]
By Lemma~\ref{lem:elliptic-grading-comparison}, the graded model of
$\Ku(X_{d,k})$ is $\mathcal C$. Lemma~\ref{lem:elliptic-product-gepner} constructs a
$\Gamma$-invariant stability condition on $\Db(E^k)$. By Lemma~\ref{lem:elliptic-equivariant-support}, it
descends to a stability condition on $\mathcal C=\Db(E^k)^\Gamma$ with total mass--Hom bounds and support on the
full numerical group. Under \eqref{eq:elliptic-coarse-category}, the grading shift satisfies
\(
V\circ(1)\simeq T_k\circ V.
\)
Hence \eqref{eq:elliptic-factor-gepner} gives
\[
 \begin{aligned}
 F\in\mathcal P_{\mathcal C}(\phi) &\Longleftrightarrow VF\in\mathcal P_A(\phi)\\
&\Longleftrightarrow T_kVF\in\mathcal P_A(\phi+2/d)\\
 &\Longleftrightarrow F(1)\in
                       \mathcal P_{\mathcal C}(\phi+2/d).
 \end{aligned}
\]
Moreover,
\[
  Z_{\mathcal C}(F(1)) =Z_A(T_kVF) =\lambda Z_A(VF) =\lambda Z_{\mathcal
 C}(F).
\]
Thus
\(
(1)\widetilde\sigma_{\mathcal C} =\widetilde\sigma_{\mathcal
C}\cdot\frac2d.
\)
Transporting $\widetilde\sigma_{\mathcal C}$ through \eqref{eq:elliptic-coarse-category} and \eqref{eq:graded-hypersurface-equivalence} gives the
required stability condition on $\Ku(X_{d,k})$.
\end{proof}

\section{Fano hypersurfaces of degrees three, four and six}
\label{sec:hypersurface-applications}

In this section, we construct numerical Gepner conditions for Fano hypersurfaces of degrees three, four
and six by passing from elliptic blocks to Fermat equations and then to smooth linear sections \cite{Orlov2009}. We also
treat quartic double fivefolds and compute the global dimensions \cite{Abuaf2020,Toda2013,KP2021}.

\subsection{Ordinary hypersurfaces}
\label{subsec:ordinary-hypersurfaces}

\begin{proposition}
\label{prop:ordinary-fermat-seeds}
For every $d\in\{3,4,6\}$ there are arbitrarily large integers $M$ for which the ordinary Fermat hypersurface
\[
 F_{d,M}=\{u_1^d+\cdots+u_M^d=0\}\subset\PP^{M-1}
\]
has a numerical stability condition on its Kuznetsov component with mass-Hom bounds, full numerical support, and satisfying
$R\widetilde\sigma=\widetilde\sigma\cdot2/d$. One may take $M=3k,2k,k$ in degrees $3,4,6$, respectively, for every $k\geq7$.
\end{proposition}

\begin{proof}
Theorem~\ref{thm:elliptic-block-seed} gives the required stability condition
for the potential $W_{d,k}$ of $k$ elliptic blocks. Its degree is $d$, and the sum of
its variable weights is $kd$. For $d=3$, all weights are one and
\[
 W_{3,k}=\sum_{j=1}^k(x_j^3+y_j^3+z_j^3).
\]
Thus the associated hypersurface is $F_{3,3k}$. For $d=4$, consider the successive coordinate restrictions
\[
W_{4,k} =\sum_{j=1}^k(x_j^4+y_j^4+z_j^2) \longmapsto W_{4,k}|_{z_1=0} \longmapsto\cdots\longmapsto
   \sum_{j=1}^k(x_j^4+y_j^4).
\]
Each removed variable has weight
\(
a=2=d/2. \) The final hypersurface is $F_{4,2k}$. For $d=6$, perform
the restrictions
\[
 \begin{aligned}
 W_{6,k}
 &=\sum_{j=1}^k(x_j^6+y_j^3+z_j^2)\\
 &\longmapsto W_{6,k}|_{z_1=\cdots=z_k=0}
   =\sum_{j=1}^k(x_j^6+y_j^3)\\
 &\longmapsto W_{6,k}|_{z_1=\cdots=z_k=y_1=\cdots=y_k=0}
   =\sum_{j=1}^k x_j^6,
 \end{aligned}
\]
removing one variable at each step. The respective cut
weights are
\(
a=3=d/2\), and \(a=2< d/2. \) The final hypersurface is $F_{6,k}$. We now verify the hypotheses of
Proposition~\ref{prop:weighted-hypersurface-restriction} at every step. Each intermediate potential has the form
\[
 W=\sum_{\nu=1}^r u_\nu^{b_\nu},
 \qquad b_\nu\geq2,
 \qquad b_\nu w_\nu=d.
\]
Hence
\[
 \left\{\frac{\partial W}{\partial u_1} =\cdots= \frac{\partial W}{\partial u_r}=0\right\} =\{0\},
\]
so its weighted projective hypersurface stack is smooth and proper. For a cut variable $q=u_\nu$, write
\[
 W=q^{b_\nu}+W',
 \qquad 0\ne W'\in\C[u_\mu\mid\mu\ne\nu].
\]
Then $\gcd(q,W)=1$, and therefore
\(
(W):q=(W).\) Thus multiplication by $q$ is injective on $P/(W)$. Moreover, the
sum of the remaining weights satisfies
\[
  \sum_\nu w_\nu\geq
 \begin{cases}
  2k>4,&d=4,\\
 k>6,&d=6,
 \end{cases}
 \qquad k\geq7.
\]
Thus all these hypersurfaces are Fano. By \eqref{eq:factorization-singularity-equivalence}
and \eqref{eq:graded-hypersurface-equivalence}, its graded singularity
category is equivalent to an admissible component of the
derived category of a smooth proper hypersurface stack. It is
therefore smooth and proper, with a finite-rank numerical Grothendieck
group.
Proposition~\ref{prop:weighted-hypersurface-restriction} now applies
successively, including the cuts of weight $d/2$. At each
step it preserves total mass--Hom bounds, full numerical support, and
the equation
\(
R\widetilde\sigma =\widetilde\sigma\cdot 2/d. \) The final
conditions lie on $\Ku(F_{3,3k})$, $\Ku(F_{4,2k})$,
and $\Ku(F_{6,k})$, respectively.
\end{proof}

\begin{lemma}
\label{lem:ordinary-waring}
Let $f\in\C[x_1,\ldots,x_m]$ be homogeneous of degree $d\geq2$. For any
sufficiently large integer $M$ in any prescribed unbounded set of positive integers, there
are linear forms $\ell_1,\ldots,\ell_M$ such that
\(
f=\sum_{j=1}^M\ell_j^d
\)
and the linear map
$\ell\colon\C^m\to\C^M$ given by these forms is injective.
\end{lemma}

\begin{proof}
By polarization,
\(
\operatorname{Sym}^d((\C^m)^\vee) = \operatorname{span}_{\C}
\{\ell^d\mid\ell\in(\C^m)^\vee\}.
\)
Thus
\[
f=\sum_{\nu=1}^r c_\nu\ell_\nu^d =\sum_{\nu=1}^r(a_\nu\ell_\nu)^d, \qquad a_\nu^d=c_\nu.
\]
Choose $\eta\in\C$ with $\eta^d=-1$ and define
\[
 L\colon\C^m\longrightarrow\C^{r+2m},
 \qquad
 x\longmapsto
 \bigl(a_1\ell_1(x),\ldots,a_r\ell_r(x),
       x_1,\eta x_1,\ldots,x_m,\eta x_m\bigr).
\]
Then
\[
\ker L=0, \qquad \sum_{\nu=1}^{r+2m}L_\nu(x)^d =f(x)+\sum_{i=1}^m(1+\eta^d)x_i^d =f(x).
\]
For any $M\geq r+2m$ in the prescribed set, replace $L$ by
\[
L\oplus0\colon\C^m\hookrightarrow\C^M.
\]
The same identity holds, proving the assertion.
\end{proof}

\begin{lemma}
\label{lem:smooth-linear-flag}
Let $F\subset\PP^{M-1}$ be a smooth hypersurface and let
$L\simeq\PP^{m-1}$ be a linear subspace such that
$X=F\cap L$ is a smooth hypersurface in $L$. There is a flag
\[
 \PP^{M-1}=L_0\supset L_1\supset\cdots\supset L_{M-m}=L
\]
with each inclusion a hyperplane inclusion and each $F\cap L_j$ smooth.
\end{lemma}

\begin{proof}
For $x\in X$, let $f_x$ be a local defining equation of $F$.
Since $X=F\cap L$ is a smooth hypersurface in $L$, every
linear subspace $\Lambda\supseteq L$ satisfies
\[
d(f_x|_\Lambda)_x|_{T_xL} =d(f_x|_L)_x\ne0.
\]
Thus $F\cap\Lambda$ is smooth along $X$.

Suppose now that $F\cap\Lambda$ is smooth and $\Lambda\supsetneq L$. The hyperplanes
$H\subset\Lambda$ containing $L$ have base locus $X$ on $F\cap\Lambda$. By Bertini theorem, a general
section $F\cap H$ is smooth away from $X$, and the preceding calculation gives smoothness along $X$.
Starting with the ambient projective space and repeating this construction gives the required
flag.
\end{proof}

\begin{proof}[Proof of Theorem~\ref{thm:intro-hypersurfaces}]
Set $m=n+2$, and let $f$ be the homogeneous equation of $X_d$. By Proposition~\ref{prop:ordinary-fermat-seeds} and
Lemma~\ref{lem:ordinary-waring}, we may choose $M$ and linear forms $\ell_1,\ldots,\ell_M$ such that
\[
f=\sum_{i=1}^M\ell_i^d, \qquad \ell=(\ell_1,\ldots,\ell_M)\colon\C^m\hookrightarrow\C^M,
\]
and $\Ku(F_{d,M})$ has the required starting condition.
Thus $\PP(\ell)$ identifies $X_d$ with a linear section
of $F_{d,M}$. By Lemma~\ref{lem:smooth-linear-flag}, there is a flag
\[
X_m=X_d\subset X_{m+1}\subset\cdots\subset X_M=F_{d,M},
\]
where $X_r\subset\PP^{r-1}$ is a smooth degree-$d$ hypersurface and $X_r$ is a hyperplane
section of $X_{r+1}$. Choose coordinates compatible with the flag, so that their equations satisfy
\(
f_r=f_{r+1}|_{x_{r+1}=0}.
\)
At each step, we have
\[
  \deg(x_{r+1})=1\leq\frac d2, \qquad f_r\ne0, \qquad
(f_{r+1}):x_{r+1}=(f_{r+1}).
\]
Hence the cut coordinate is regular in the hypersurface ring. Moreover, $r\geq
m>d$, so every $X_r$ is Fano. By \eqref{eq:factorization-singularity-equivalence} and
\eqref{eq:graded-hypersurface-equivalence}, the categories $\MF^\Z(f_r)$ are smooth and proper, with
finite-rank numerical Grothendieck groups. Applying Proposition~\ref{prop:weighted-hypersurface-restriction}
successively gives a stability condition $\widetilde\sigma$ on $\MF^\Z(f)$ satisfying
\(
R\widetilde\sigma=\widetilde\sigma\cdot2/d,
\)
with total mass--Hom bounds and full numerical support.
Transport through \eqref{eq:graded-hypersurface-equivalence}
gives the required condition $\widetilde\sigma_X$ on $\Ku(X_d)$.

Finally, the Serre functor is
\[
\Sfun \simeq R^{d-m}[m-2] \simeq R^{-m}[m]
\]
by \cite{Toda2013,KP2021}. The grading equation therefore gives
\[
\Sfun\mathcal P(\phi)=\mathcal P(\phi+s), \qquad s=m-\frac{2m}{d}
 =\frac{(n+2)(d-2)}{d}.
\]
Lemma~\ref{lem:serre-phase-gldim} now yields
\[
\gldim(\widetilde\sigma_X) =\frac{(n+2)(d-2)}{d}.
\]
\end{proof}

\subsection{Quartic double fivefolds}
\label{subsec:quartic-double-fivefolds}

The weighted construction also applies when one retains a quadratic variable. Let
\[
  Y=\{z^2+f_4(x_1,\ldots,x_6)=0\} \subset\PP(1,1,1,1,1,1,2),
\]
where the branch quartic $\{f_4=0\}\subset\PP^5$ is smooth. The hypersurface avoids the point with nontrivial stabilizer in
the ambient weighted projective space. Thus $Y$ is a smooth projective variety. We set
\[
\Ku(Y)=\langle\OO_Y,\OO_Y(1),\OO_Y(2),\OO_Y(3)\rangle^\perp.
\]
The semiorthogonal decomposition and its Calabi--Yau structure
are recalled in \cite[Section~3.2]{Abuaf2020}.

\begin{lemma}
\label{lem:quartic-double-fermat-seeds}
For every sufficiently large integer $k$, the weighted hypersurface
\[
 Y_{2k}=\{z^2+u_1^4+\cdots+u_{2k}^4=0\} \subset\PP(\underbrace{1,\ldots,1}_{2k},2)
\]
has a numerical stability condition on $\Ku(Y_{2k})$ with mass-Hom bounds, full
numerical support, and
$R_{Y_{2k}}\widetilde\sigma=\widetilde\sigma\cdot1/2$.
\end{lemma}

\begin{proof}
Theorem~\ref{thm:elliptic-block-seed} gives the required stability condition for
\[
W_k=\sum_{j=1}^k(x_j^4+y_j^4+z_j^2).
\]
Consider the successive restrictions
\[
 W_r=\sum_{j=1}^k(x_j^4+y_j^4)+\sum_{j=1}^r z_j^2, \qquad W_{r-1}=W_r|_{z_r=0}, \qquad r=k,\ldots,2.
\]
Each cut has weight $2=4/2$, so it lies at the endpoint of
Proposition~\ref{prop:weighted-hypersurface-restriction}. Every $W_r$
satisfies
\[
 \operatorname{Crit}(W_r)=\{0\}, \qquad \sum_\nu w_\nu=2k+2r\geq2k+2>4.
\]
Its weighted projective hypersurface is therefore a smooth proper Fano stack. Moreover,
\[
 W_r|_{z_r=0}=W_{r-1}\ne0,
 \qquad
 (W_r):z_r=(W_r),
\]
so $z_r$ is regular in the hypersurface ring. By \eqref{eq:factorization-singularity-equivalence} and
\eqref{eq:graded-hypersurface-equivalence}, the corresponding graded singularity categories are smooth and proper.
Proposition~\ref{prop:weighted-hypersurface-restriction} thus applies at every step, preserving total mass--Hom bounds, full
numerical support, and the grading equation
\(
R\widetilde\sigma=\widetilde\sigma\cdot1/2.
\)
The final potential is
\[
W_1=G_0(u_1,\ldots,u_{2k})+z^2, \qquad G_0=\sum_{j=1}^{2k}u_j^4.
\]
Transporting the resulting stability condition through \eqref{eq:graded-hypersurface-equivalence} gives the required
stability condition on $\Ku(Y_{2k})$.
\end{proof}

\begin{lemma}
\label{lem:quartic-double-smooth-flag}
Every smooth quartic double fivefold $Y$ is obtained from $Y_{2k}$, for some sufficiently
large $k$, by a sequence of coordinate cuts of weight one. Every
intermediate hypersurface is smooth and Fano, and every
cut coordinate is regular in its hypersurface ring.
\end{lemma}

\begin{proof}
By Lemma~\ref{lem:ordinary-waring}, after increasing $k$ we may write
\[
f_4=\sum_{i=1}^{2k}\ell_i^4, \qquad (\ell_1,\ldots,\ell_{2k})\colon\C^6\hookrightarrow\C^{2k}.
\]
Thus $\{f_4=0\}$ is a linear section of the Fermat quartic $\{G_0=0\}$.
Lemma~\ref{lem:smooth-linear-flag} gives a flag of smooth
quartic hypersurfaces
\[
X_6\subset X_7\subset\cdots\subset X_{2k}, \qquad X_6=\{f_4=0\}, \quad X_{2k}=\{G_0=0\}.
\]
Choose coordinates compatible with the flag and write
\[
  X_r=\{G_r=0\}\subset\PP^{r-1}, \qquad G_{r-1}=G_r|_{u_r=0}.
\]
Retaining the variable $z$ of weight two, set
\[
 W_r=G_r+z^2, \qquad Y_r=\{W_r=0\}\subset\PP(1^r,2), \qquad 6\leq r\leq2k.
\]
Since $X_r$ is smooth,
\[
 \operatorname{Crit}(W_r)
 =\{(u,z)\mid \nabla G_r(u)=0,\ 2z=0\}
 =\{0\}.
\]
Moreover, $Y_r$ avoids the stacky point $[0:\cdots:0:1]$. Projection onto the
ordinary coordinates therefore realizes $Y_r$ as a smooth double cover of
$\PP^{r-1}$ branched along $X_r$. Adjunction then gives
\[
\omega_{Y_r}\simeq\mathcal O_{Y_r}(4-r-2) =\mathcal O_{Y_r}(2-r),
\]
so every $Y_r$ is Fano. The successive cuts satisfy
\[
 W_r|_{u_r=0}=W_{r-1}\ne0,
 \qquad
 \deg(u_r)=1\leq\frac42.
\]
Hence $\gcd(u_r,W_r)=1$, and multiplication by $u_r$ is injective on the
corresponding hypersurface ring. Finally,
\eqref{eq:factorization-singularity-equivalence} and
\eqref{eq:graded-hypersurface-equivalence} identify the intermediate graded
singularity categories with admissible components of $\Db(Y_r)$. They are therefore smooth
and proper, with finite-rank numerical Grothendieck groups.
\end{proof}

\begin{proposition}
\label{prop:quartic-double-existence}
Every smooth quartic double fivefold $Y$ admits a numerical stability condition $\widetilde\sigma_Y$ on $\Ku(Y)$ with mass-Hom bounds, full
numerical support, and
\[
 R_Y\widetilde\sigma_Y=\widetilde\sigma_Y\cdot\frac12.
\]
\end{proposition}

\begin{proof}
By Lemma~\ref{lem:quartic-double-smooth-flag}, choose a flag
\[
Y=Y_6\subset Y_7\subset\cdots\subset Y_{2k}.
\]
Lemma~\ref{lem:quartic-double-fermat-seeds} gives a stability
condition on $\Ku(Y_{2k})$ with the stated properties.
Each step is a coordinate cut with
\(
a=1<d/2=2.
\)
Applying Proposition~\ref{prop:weighted-hypersurface-restriction}
successively therefore gives a stability condition $\widetilde\sigma_Y$ on
$\Ku(Y)$ with total mass--Hom bounds, full numerical support, and
\(
R_Y\widetilde\sigma_Y =\widetilde\sigma_Y\cdot1/2.
\)
\end{proof}

\section{Complete intersections of three quadrics}
\label{sec:three-quadrics}

In this section, we construct stability conditions for smooth Fano intersections of three
quadrics using \cite{KuzQuadrics2008,CO2026}. We verify the geometric hypotheses for twisted
surfaces in even dimensions and use a Clifford-module surface tilt in odd
dimensions, establishing full numerical support in both cases
\cite{Koseki2023,MacriSchmidt2017}.

\subsection{Even-dimensional intersections}

\begin{theorem}\label{thm:even-three-quadrics}
Let $m\geq2$ and let
\(
X\subset\Pp^{2m+3}_{\C}
\)
be a smooth complete intersection of three quadrics. Set
\[
 \Ku(X)=
 \langle\mathcal O_X,\mathcal O_X(1),\ldots,
                 \mathcal O_X(2m-3)\rangle^\perp.
\]
There are a smooth projective surface $S$, a Brauer class $\alpha\in\operatorname{Br}(S)$, and an exact equivalence
\[
 \Phi\colon\Db(S,\alpha)\xrightarrow{\ \sim\ }\Ku(X).
\]
In particular, $\Ku(X)$ admits a numerical Bridgeland stability condition $\widetilde\sigma_X$ with the
support property on its full numerical Grothendieck group.
\end{theorem}

This theorem follows from \cite[Proposition~8.7(ii)]{CO2026}
and the construction of stability conditions on twisted surfaces. The argument below proves that its assumptions hold for every smooth complex
complete intersection in the statement.

\begin{lemma}\label{lem:net-geometric-checks}
Let $V$ be a complex vector space, and let $q_0,q_1,q_2$
be linearly independent quadratic forms defining a smooth
complete intersection
\[
 X=\{q_0=q_1=q_2=0\}\subset\Pp(V)
\]
of codimension three. Then the universal quadric
\[
 Q=\left\{(x,[t_0:t_1:t_2])\in\Pp(V)\times\Pp^2
 \ \middle|\ \sum_{i=0}^2t_iq_i(x)=0\right\}
\]
is smooth, and the projection $\rho\colon Q\to\Pp^2$ is a flat quadric bundle
whose fibres have corank at most two. The locus in $\Pp^2$ where the fibres
have corank two is a finite reduced scheme.
\end{lemma}

\begin{proof}
Since $q_0,q_1,q_2$ are linearly independent, every fibre is defined by a
nonzero quadratic form and has the same Hilbert polynomial. Thus
$\rho$ is a flat quadric bundle. On the chart $t_0\ne0$,
write
\(
q_s=q_0+s_1q_1+s_2q_2,\) for \(s_i=t_i/t_0.
\)
If $(x,s)\in Q$ were singular, then
\[
  q_1(x)=q_2(x)=0, \qquad q_s(x)=0, \qquad d_xq_s=0.
\]
Hence $x\in X$ and
\[
 d_xq_0+s_1d_xq_1+s_2d_xq_2=0.
\]
After a local trivialization of $\mathcal O_{\Pp(V)}(2)$, smoothness and
codimension three of $X$ give
\[
\dim\operatorname{span}\{d_xq_0,d_xq_1,d_xq_2\}=3,
\]
which is a contradiction. The same argument on the other charts
proves that $Q$ is smooth. Fix $t\in\Pp^2$, and let $K=\ker(q_t)$ denote the kernel
of the associated symmetric bilinear form.
Choose local coordinates $s_1,s_2$ centered at $t$, and set
\[
 b_i=\left.\frac{\partial q_s}{\partial s_i}\right|_{s=0,K} \in\operatorname{Sym}^2K^\vee.
\]
The fibre differential vanishes at every point of $\Pp(K)$. Smoothness of $Q$ therefore implies
\[
  \{[x]\in\Pp(K)\mid b_1(x)=b_2(x)=0\}=\varnothing.
\]
If $\dim K\geq3$, two quadratic forms on $\Pp(K)$ necessarily have a
common zero. Thus
\(
\dim K\leq2 \) for every fibre. Suppose $\dim K=2$. Choose a complement of $K$
and write the local symmetric matrix of the family as
\[
A(s)=
 \begin{pmatrix}
  U(s)&B(s)\\
B(s)^{\mathsf t}&D(s)
 \end{pmatrix},
 \]
where $U(0)$ is invertible and $B(0)=D(0)=0$. After shrinking the neighbourhood, $U(s)$ is invertible, and
$A(s)$ is congruent to
\[
U(s)\oplus C(s), \qquad C(s)=D(s)-B(s)^{\mathsf t}U(s)^{-1}B(s).
\]
Hence the corank-two locus $Z$ is defined locally by the three entries of the symmetric
matrix $C(s)$, and
\[
  T_tZ =\ker\left( dC_t\colon T_t\Pp^2\longrightarrow\operatorname{Sym}^2K^\vee \right).
\]
Since $B(0)=0$, we have
\[
dC_t\left(\frac{\partial}{\partial s_i}\right)=b_i.
\]
The forms $b_1,b_2$ have no common zero on $\Pp(K)\simeq\Pp^1$. They are therefore linearly independent,
since every nonzero binary quadratic has a projective zero. Consequently $dC_t$ is injective.

Let $I\subset\mathcal O_{\Pp^2,t}$ be the ideal generated by the entries of $C$, and let $\mathfrak m_t$ be the maximal ideal. Injectivity of
$dC_t$ gives
\(
I+\mathfrak m_t^2=\mathfrak m_t. \) By Nakayama's lemma, $I=\mathfrak m_t$. Thus every point of $Z$ is reduced and isolated. Since $Z$ is closed in
$\Pp^2$, it is a finite reduced scheme.
\end{proof}

\begin{lemma}\label{lem:twisted-surface-full-support}
Let $S$ be a smooth projective complex surface and let
$\alpha\in\operatorname{Br}(S)$. The category
$\Db(S,\alpha)$ admits a numerical Bridgeland stability
condition with the support property on its full
numerical Grothendieck group.
\end{lemma}

\begin{proof}
Represent $\alpha$ by a finite banded gerbe over $S$, and identify
$\operatorname{Coh}(S,\alpha)$ with its category of
weight-one coherent sheaves. Write the normalized twisted
Chern character as
\[
 v(E)=(r(E),c(E),s(E)) \in\Q\oplus\operatorname{NS}(S)_\Q\oplus\Q, \qquad \Delta(E)=c(E)^2-2r(E)s(E).
\]
The normalization uses the inverse of pullback on rational
Chow groups. In particular, a weight-one skyscraper sheaf
has character $(0,0,1)$. For twisted perfect complexes $E,F$, the weight-zero
complex $E^\vee\otimes F$ descends to $S$. Riemann--Roch
therefore gives
\[
\chi(E,F) =\int_S v(E)^\vee v(F)\operatorname{td}(S), \qquad (r,c,s)^\vee=(r,-c,s).
\]
Moreover, the image of $v$ spans
\(
V=\R\oplus\operatorname{NS}(S)_\R\oplus\R. \) Indeed, choose a locally free twisted sheaf $G$ of rank $r_G>0$, which exists
because $\alpha$ is a Brauer class. For every integral divisor $D$,
\[
 v(G(D))
 =\left(r_G,c_G+r_GD,
 s_G+c_G\cdot D+\frac{r_G}{2}D^2\right),
\]
and hence
\[
 v(G(D))-v(G(-D))
 =(0,2r_GD,2c_G\cdot D).
\]
Together with $(0,0,1)$ and $v(G)$, these classes span $V$. The surface Euler
pairing on $V$ is nondegenerate by nondegeneracy of the divisor intersection form.
Thus
\(
\ker v=\operatorname{rad}\chi,\) and \(\Knum(S,\alpha)_\R\simeq V.
\)
Choose twisted objects $G_1,\ldots,G_t$ whose characters form a basis of $\Q\oplus\operatorname{NS}(S)_\Q\oplus\Q$. The map
\[
V\longrightarrow\R^t,\qquad w\longmapsto\bigl(\chi(v(G_i),w)\bigr)_{i=1}^t
\]
is an isomorphism with rational coefficients. Since $\chi(G_i,E)\in\Z$, the image of $v$ is
contained in a fixed lattice in $V$. It spans $V$, so it is itself a
lattice. These constructions agree with \cite[Section~2.2.7]{Lieblich2007}. It therefore
suffices to prove support on $V$. Choose an integral ample divisor $H$, and set
\(
h=H^2,\) and
 \(
\delta(r,c,s)=H\cdot c. \) For a torsion-free twisted sheaf $F$, set
\[
\mu_H(F)=\frac{H\cdot c(F)}{hr(F)},
\]
and write $\mu_H^+(F)$ and $\mu_H^-(F)$ for its maximal and minimal HN slopes. Using the
slope HN filtration, define
\[
 \begin{aligned}
\mathcal F &=\{F\text{ torsion-free}\mid\mu_H^+(F)\leq0\},\\
\mathcal T &=\{E\mid E\text{ is torsion, or } \mu_H^-(E/E_{\mathrm{tor}})>0\}.
 \end{aligned}
\]
Let
\[
\mathcal A=\langle\mathcal F[1],\mathcal T\rangle, \qquad Z_a(r,c,s)=-s+\frac{a^2h}{2}r+ia\delta(r,c,s), \qquad a\in\R_{>0}.
\]

For a slope-semistable torsion-free twisted sheaf $F$, the twisted Bogomolov inequality gives
\(
\Delta(F)\geq0 \) by \cite[Proposition~3.2.3.13]{Lieblich2007},
whose correction term vanishes in characteristic zero. Together with the
Hodge index theorem, this shows that $Z_a$ is a stability function on $\mathcal A$.
Indeed, $\delta(E)\geq0$ for $E\in\mathcal A$. If $\delta(E)=0$, its cohomology sheaves
are a slope-semistable torsion-free sheaf $F$ of slope zero and a zero-dimensional sheaf $T$. In
this case
\[
s(F)\leq\frac{c(F)^2}{2r(F)}\leq0,
\]
and
\[
 Z_a(F[1])=s(F)-\frac{a^2h}{2}r(F)<0, \qquad Z_a(T)=-s(T)<0
\]
whenever the corresponding sheaf is nonzero.

The category $\operatorname{Coh}(S,\alpha)$ is noetherian and has slope HN filtrations. Torsion-free
twisted sheaves have locally free double duals with zero-dimensional quotient. Together with discreteness of
$\delta$, these are the inputs in the proof of \cite[Lemma~6.17]{MacriSchmidt2017}, which therefore shows that $\mathcal
A$ is noetherian. By \cite[Proposition~4.10]{MacriSchmidt2017}, $Z_a$ has the HN property for every
real $a>0$. This is the surface tilt construction of \cite[Section~2]{ArcaraBertram2013}.

We now construct a support form on the full space $V$. There is a constant
$C_H\geq0$ such that
\(
D^2+C_H(H\cdot D)^2\geq0\) for every effective rational divisor $D$. To see
this, choose a basis $D_1,\ldots,D_\rho$ of $\operatorname{NS}(S)_\R$
and $\varepsilon>0$ such that $H\pm\varepsilon D_j$ are ample. For effective
$D$,
\[
(H\pm\varepsilon D_j)\cdot D\geq0, \qquad |D_j\cdot D|\leq\varepsilon^{-1}H\cdot D.
\]
Nondegeneracy of the intersection pairing then gives
$\|D\|\leq C(H\cdot D)$ for any fixed norm. The asserted
inequality follows by continuity of the intersection form. Fix such a constant and define
\[
Q_H(r,c,s)=c^2-2rs+C_H\delta(r,c,s)^2.
\]
For $v\in\ker Z_a$,
\(
H\cdot c=0,\) and \( s={a^2h}r/{2},
\)
so
\(
Q_H(v)=c^2-a^2hr^2. \) The Hodge index theorem therefore gives
\(
Q_H(v)<0,\) for \( 0\ne v\in\ker Z_a.
\) 

It remains to prove $Q_H(v(E))\geq0$ for $Z_a$-semistable objects. We follow the discrete degree induction of
\cite[Theorem~6.13 and Section~6.3]{MacriSchmidt2017}, recording the sheaf inequalities needed in the
twisted case. For a torsion twisted sheaf $T$, the class $c(T)$ is effective and $r(T)=0$, so
\[
 Q_H(v(T))
 =c(T)^2+C_H(H\cdot c(T))^2\geq0.
\]
For a slope-semistable torsion-free twisted sheaf $F$,
\[
Q_H(v(F)) =\Delta(F)+C_H(H\cdot c(F))^2\geq0.
\]
If $E\in\mathcal A$ has $\delta(E)=0$, write
\[
  F=H^{-1}(E), \qquad T=H^0(E), \qquad \ell=s(T)\geq0.
\]
Here $F$ is slope-semistable of slope zero and $T$ is
zero-dimensional. Since $v(E)=v(T)-v(F)$,
\[
 Q_H(v(E))=\Delta(F)+2r(F)\ell\geq0.
\]

Suppose $\delta(E)>0$, and assume the assertion for smaller positive degrees. Negative definiteness on $\ker Z_a$ implies
that, for every phase $\phi$, the cone
\[
 \left\{
 v\in V\ \middle|\
 Z_a(v)\in\R_{>0}e^{i\pi\phi},\quad Q_H(v)\geq0
 \right\}
\]
is convex. Splitting off $\ker Z_a$ reduces its defining
inequality to a Lorentz-cone inequality. If $E$ is strictly
semistable, there is an exact sequence in $\mathcal A$
\[
 0\longrightarrow F\longrightarrow E\longrightarrow G \longrightarrow0
\]
whose nonzero terms have the same phase. Degree-zero objects have phase one, so
\[
0<\delta(F),\delta(G)<\delta(E).
\]
Induction and convexity give $Q_H(v(E))\geq0$.

It therefore remains to exclude a stable object $E$ with $Q_H(v(E))<0$. For the
wall argument use
\[
\widehat Z_a=-s+\frac{a^2h}{2}r+i\delta.
\]
It has the same semistable objects and kernel as $Z_a$. The local norm and mass estimates in the proof of
Proposition~\ref{prop:nc-surface-stability} apply with $Q_H$ and $V$. They use discreteness of the degree and
negative definiteness on the central-charge kernel. The induction hypothesis supplies the norm estimate for the HN
factors of every potential destabilizing subobject. Thus only finitely many numerical classes can destabilize near any
parameter, and each equality of phases is linear in $a^2$. These estimates are uniform on compact intervals in $\R_{>0}$, so a
finite endpoint of an interval of stability is an actual wall. At such a wall, the same-phase exact sequence
above and induction would give $Q_H(v(E))\geq0$. Hence $E$ remains semistable for every
larger $a$.

By \cite[Lemma~6.18]{MacriSchmidt2017}, there is only one of the
following possibilities
\begin{enumerate}
 \item $E$ is a slope-semistable torsion-free sheaf;
 \item $E$ is a torsion sheaf;
 \item $H^{-1}(E)=F$ is slope-semistable torsion-free and
$H^0(E)=T$ is zero-dimensional.
\end{enumerate}
The comparison uses the cohomology triangle and the slope filtration, so it applies to twisted sheaves as
well. The first two cases were checked above. In the third,
\[
Q_H(v(E)) =\Delta(F)+2r(F)s(T)+C_H(H\cdot c(F))^2\geq0.
\]
Each case contradicts $Q_H(v(E))<0$. This completes the induction. As in the proof of
\cite[Theorem~6.13]{MacriSchmidt2017}, the first-wall argument uses the induction
hypothesis for possible destabilizing classes and does not
assume full support. Thus
\(
Q_H(v(E))\geq0\) {for every semistable }$E$, and $Q_H|_{\ker Z_a}<0$, which proves the support property on $V$, including all divisor directions in
$H^\perp$, and hence on the full numerical lattice. Therefore $(\mathcal A,Z_a)$ defines the required numerical stability condition.
\end{proof}

\begin{proof}[Proof of Theorem~\ref{thm:even-three-quadrics}]
Choose quadratic equations $q_0,q_1,q_2$ for $X$, and let
\(
\rho\colon Q\longrightarrow\Pp^2
\)
be their universal quadric bundle.
By Lemma~\ref{lem:net-geometric-checks}, $Q$ is smooth,
every fibre has corank at most two, and the corank-two
locus is finite and reduced. Set $\ell=m+1$. Then
\[
X\subset\Pp^{2\ell+1}, \qquad \dim Q_t=2\ell, \qquad 3\leq\ell+1.
\]
Thus \cite[Proposition~8.7(ii)]{CO2026} applies and gives a smooth projective surface $S$, a
Brauer class $\alpha$, and a semiorthogonal decomposition
\[
\Db(X)= \left\langle \Db(S,\alpha), \mathcal O_X(1),\ldots,\mathcal
   O_X(2\ell-4) \right\rangle.
\]
Since $2\ell-4=2m-2$, tensoring by $\mathcal O_X(-1)$
identifies the first component with $\Ku(X)$ and yields
an equivalence
\(
\Phi\colon\Db(S,\alpha)\xrightarrow{\sim}\Ku(X). \) By Lemma~\ref{lem:twisted-surface-full-support}, choose
\(
\widetilde\sigma_S=(\mathcal P_S,Z_S)
\in\Stab_{\mathrm{num}}\bigl(\Db(S,\alpha)\bigr)
\)
with support on the full numerical group. Define
\[
\mathcal P_X(\phi)=\Phi\bigl(\mathcal P_S(\phi)\bigr), \qquad
Z_X=Z_S\circ\Phi_*^{-1}.
\]
The equivalence preserves shifts, Hom groups, and distinguished triangles, so
it transports the slicing axioms and HN filtrations. Since $\Phi$ preserves the Euler
pairing, it induces
\[
  \Phi_{\mathrm{num}}\colon \Knum\bigl(\Db(S,\alpha)\bigr)
\xrightarrow{\sim}\Knum\bigl(\Ku(X)\bigr).
\]
For a support form $Q_S$ of $\widetilde\sigma_S$, set
\(
Q_X=Q_S\circ\Phi_{\mathrm{num}}^{-1}.
\)
Then
\[
\ker Z_X=\Phi_{\mathrm{num}}(\ker Z_S), \qquad Q_X(v_X(E))
=Q_S\bigl(v_S(\Phi^{-1}E)\bigr).
\]
Hence $Q_X$ is negative definite on $\ker Z_X$ and
nonnegative on every $\widetilde\sigma_X$-semistable
class. Thus
\[
 \widetilde\sigma_X=(\mathcal P_X,Z_X) \in\Stab_{\mathrm{num}}\bigl(\Ku(X)\bigr).
\]
\end{proof}

\subsection{A surface construction for sheaves of algebras}
\label{subsec:nc-surface-tilt}

The Clifford category associated with a net of quadrics is a derived category of modules over an algebra on a surface. To show the existence of
stability conditions, a sufficiently weak Bogomolov inequality is enough. We give the argument in the rank-three case, which is the case needed below.

\begin{lemma}\label{lem:nc-slope-spread}
Let $\mathcal B$ be a coherent unital, $\mathcal O_{\mathbb P^2}$-central algebra which is
locally free as an $\mathcal O_{\mathbb P^2}$-module. Choose an integer $N\geq 0$
such that $\mathcal B(N)$ is globally generated. Every slope semistable
torsion-free $\mathcal B$-module $E$ satisfies
\begin{equation}\label{eq:nc-weak-bg}
 \Delta(E):=c_1(E)^2-2\operatorname{rk}(E)\operatorname{ch}_2(E)
 \ \geq -N^2\operatorname{rk}(E)^2.
\end{equation}
\end{lemma}

\begin{proof}

Let $H$ be the hyperplane class on $\mathbb P^2$, so
$H^2=1$, and write $\mu(F)=H\cdot c_1(F)/r(F)$.
Let $F\subset E$ be the first ordinary
HN factor, and set
\(
G=\operatorname{im}(\mathcal B\otimes F\to E).
\)
Then $G$ is a nonzero $\mathcal B$-submodule. By the
choice of $N$, there is a surjection
\(
F(-N)^{\oplus b}\twoheadrightarrow G.
\)
Since $F$ is slope semistable,
\[
\mu_{\max}(E)-N =\mu(F)-N \leq\mu_{\min}(G) \leq\mu(G) \leq\mu(E).
\]
Thus $\mu_{\max}(E)\leq\mu(E)+N$. Set $A=E^{\vee\vee}$. The $\mathcal B$-action extends uniquely to $A$ by local freeness of $\mathcal B$ and
extension of morphisms between reflexive sheaves across a finite subset. Moreover, $A$ is slope semistable as a $\mathcal B$-module.
Indeed, for every nonzero $\mathcal B$-submodule $M\subset A$, the quotient $M/(M\cap E)$ is zero-dimensional, so
\[
 \mu(M)=\mu(M\cap E)\leq\mu(E)=\mu(A).
\]

The dual $A^\vee$ is a right $\mathcal B$-module and is slope semistable. To see this, let $G\subset A^\vee$ be a saturated
right submodule and set
\(
Q=\operatorname{im}(A\to G^\vee).
\)
The cokernel of $Q\subset G^\vee$ is zero-dimensional. Since $Q$ is a
left-module quotient of $A$,
\(
-\mu(G)=\mu(Q)\geq\mu(A).
\)
Applying the maximal-slope estimate to $A$ and $A^\vee$ therefore gives
\[
  \mu(E)-N \leq\mu_{\min}(A) \leq\mu_{\max}(A) \leq\mu(E)+N.
\]

Let $E_1,\ldots,E_s$ be the ordinary HN factors of $A$, and write $r=r(A)$, $r_i=r(E_i)$,
$w_i=r_i/r$ and $\mu_i=\mu(E_i)$. Since $H^2=1$ and $c_1(E_i)=r_i\mu_iH$, additivity of the
Chern character and the ordinary Bogomolov inequality give
\[
 \begin{aligned}
\frac{\Delta(A)}{r^2} &=\sum_iw_i\frac{\Delta(E_i)}{r_i^2}
-\sum_iw_i\bigl(\mu_i-\mu(E)\bigr)^2\\
&\geq-N^2.
 \end{aligned}
\]
Finally, if $l=\operatorname{length}(A/E)$, then
\[
  \Delta(E)=\Delta(A)+2rl\geq-N^2r^2.
\]
This proves the assertion.
\end{proof}

Let $\operatorname{Coh}(\mathbb P^2,\mathcal B)$ denote the category
of coherent left $\mathcal B$-modules. This category is noetherian, and the
torsion subsheaf of any module is a submodule. Then HN filtrations with respect to $\mu$
can be constructed directly. For a fixed torsion-free module $E$, its submodules are ordinary
subsheaves, so their degrees are bounded above and their ranks belong to a finite set. Since degrees are
discrete, there is a submodule of maximal slope and, among these, maximal rank. Its saturation is again a
$\mathcal B$-submodule and has slope at least that of the chosen submodule. Maximality of the slope forces equality, so we
may choose the submodule saturated. Passing to the torsion-free quotient and repeating gives a finite filtration.
Its factors are semistable, and their slopes strictly decrease by maximality of the chosen rank.

For $\beta\in\mathbb Q$, let $(\mathcal T_\beta,\mathcal F_\beta)$
be the slope torsion pair, with torsion modules in
$\mathcal T_\beta$ and torsion-free slope factors of
$\mathcal F_\beta$ having slope at most $\beta$. Set
\(
 \mathcal H_\beta
 =\langle\mathcal F_\beta[1],\mathcal T_\beta\rangle.
\)
For an object $E$, write
\[
 r(E)=\operatorname{ch}_0(E),\qquad c(E)=H\operatorname{ch}_1(E),\qquad s(E)=\operatorname{ch}_2(E),
\]
where $H^2=1$. Use $c^\beta=c-\beta r$ and $s^\beta=s-\beta c+\beta^2r/2$.

\begin{proposition}\label{prop:nc-surface-stability}
Under the assumptions of Lemma~\ref{lem:nc-slope-spread}, let $a>0$ with
$a^2\in\mathbb Q$ and let $\beta\in\mathbb Q$. The pair
\begin{equation}\label{eq:nc-surface-charge}
\widetilde\sigma_{a,\beta} =\left(\mathcal H_\beta,Z_{a,\beta}\right),\qquad
 Z_{a,\beta}=-s^\beta+\frac{a^2+N^2}{2}r +\sqrt{-1}\,c^\beta,
\end{equation}
is a stability condition with support property on the lattice given by the image of $(r,c,s)$. Every $\widetilde\sigma_{a,\beta}$-semistable object
satisfies
\begin{equation}\label{eq:nc-tilt-bg}
 Q(E):=c(E)^2-2r(E)s(E)+N^2r(E)^2\geq0.
\end{equation}
If $(r,c,s)$ identifies the rational numerical Grothendieck group with $\mathbb
Q^3$, this stability condition has full numerical support.
\end{proposition}

\begin{proof}
We first verify the categorical inputs for the surface tilt
construction of \cite[Theorem~6.6]{Koseki2023} and the
rational degree induction of
\cite[Section~6.3]{MacriSchmidt2017}. For
$E\in\mathcal H_\beta$, we have $c^\beta(E)\geq0$. If
$c^\beta(E)=0$, then
\(
F=H^{-1}(E)\) {is slope semistable of slope }\(\beta\), and
\(T=H^0(E)\) {is zero-dimensional}. Lemma~\ref{lem:nc-slope-spread} gives
\[
 s^\beta(F)\leq\frac{N^2}{2}r(F), \qquad Z_{a,\beta}(F[1])\leq-\frac{a^2}{2}r(F), \qquad
Z_{a,\beta}(T)=-\operatorname{length}(T).
\]
Consequently,
\[
  Z_{a,\beta}(E)\in \{z\in\C\mid\operatorname{Im}z>0\}\cup\R_{<0} \text{ for
 }0\ne E\in\mathcal H_\beta.
\]

We next prove that $\mathcal H_\beta$ is noetherian.
Consider an increasing chain
\[
A_1\subset A_2\subset\cdots\subset M.
\]
The values $c^\beta(A_j)$ are discrete, nonnegative, and bounded by
$c^\beta(M)$, so they eventually stabilize. Quotienting by a fixed
term, we may assume $c^\beta(A_j)=0$. Set $Q_j=M/A_j$. We
have
\(
H^{-1}(A_j)\subset H^{-1}(M)
\)
and the kernels of
\(
H^0(M)\twoheadrightarrow H^0(Q_j)
\)
eventually stabilize. The cohomology sequences then give fixed modules $U,V$ and exact
sequences
\[
  0\longrightarrow U\longrightarrow H^{-1}(Q_j) \longrightarrow H^0(A_j)\longrightarrow
 V\longrightarrow0.
\]
Set
\(
 K_j=\ker\bigl(H^0(A_j)\to V\bigr).
\)
Since $H^0(A_j)$ is zero-dimensional, there is the following short exact
sequence
\[
 0\longrightarrow U\longrightarrow H^{-1}(Q_j) \longrightarrow K_j\longrightarrow0.
\]
The middle term is torsion-free and agrees with $U$ away from finitely many
points. Hence,
\[
U\subset H^{-1}(Q_j)\subset U^{\vee\vee}, \qquad \operatorname{length}(K_j)
\leq\operatorname{length}(U^{\vee\vee}/U).
\]

For $C_j=A_{j+1}/A_j$, stabilization of $H^{-1}(A_j)$ makes $H^{-1}(C_j)$ a submodule of the zero-dimensional module
$H^0(A_j)$. Since $H^{-1}(C_j)$ is torsion-free, it vanishes. Thus $C_j$ is zero-dimensional, and
\[
  \operatorname{length}(K_{j+1}) -\operatorname{length}(K_j) =\operatorname{length}(C_j).
\]
The length bound forces the chain to stabilize.

Noetherianity and discreteness of $c^\beta$ give the
HN property by the maximal-destabilizing
subobject criterion by
\cite[Proposition~4.10]{MacriSchmidt2017}.
This applies to every real $a>0$ and requires no
discreteness of the real part of the charge.

It remains to prove \eqref{eq:nc-tilt-bg}. Consider the additive coordinate
\(
 s'=s-{N^2r}/{2}.
\)
Then
\[
Q=(c^\beta)^2-2r(s')^\beta, \qquad
   Z_{a,\beta}=-(s')^\beta+\frac{a^2}{2}r+ic^\beta.
\]
In particular,
\[
v\in\ker Z_{a,\beta} \Longrightarrow c^\beta(v)=0,\quad (s')^\beta(v)=\frac{a^2}{2}r(v), \quad Q(v)=-a^2r(v)^2.
\]
Thus $Q$ is negative definite on the one-dimensional kernel for every $a>0$. We argue by induction on the
discrete values of $c^\beta(E)$. If $c^\beta(E)=0$, write $F=H^{-1}(E)$ and $T=H^0(E)$ as above. Then
\[
  Q(E) =Q(F)+2r(F)\operatorname{length}(T)\geq0
\]
by Lemma~\ref{lem:nc-slope-spread}. For positive degree, every genuine wall has semistable
factors of strictly smaller positive degree. Indeed,
degree-zero objects have phase one, so they cannot occur
as same-phase factors of an object of positive degree. Fix $a_0>0$, and suppose $E$ is semistable at $a_0$.
Write
\[
  d=c^\beta(E)>0, \qquad k=\frac{\operatorname{Re}Z_{a_0,\beta}(E)}d, \qquad
  B=1+2\max\{0,-k\}.
\]
A degree-zero subobject would contradict semistability. A proper subobject of degree $d$ has a degree-zero quotient of phase one and cannot
destabilize at any parameter. It therefore suffices to consider
\(
F\subset E,\) with \(0<c^\beta(F)<d. \) Let $G_1,\ldots,G_t$ be its HN factors at $a_0$. Their degrees are positive and
smaller than $d$, and their phases are at most $\phi_{a_0,\beta}(E)$. Hence
\(
\operatorname{Re}Z_{a_0,\beta}(G_j) \geq k\,c^\beta(G_j).
\)
Using $|x+iy|\leq|x|+y$ for $y\geq0$, we obtain
\[
 m_{a_0,\beta}(F)
 \leq\operatorname{Re}Z_{a_0,\beta}(F)+Bd.
\]
By induction, $Q(G_j)\geq0$. Negative definiteness on
$\ker Z_{a_0,\beta}$ therefore gives a constant $C$ with
\[
\sum_j\|v(G_j)\| \leq C\,m_{a_0,\beta}(F).
\]
Since the charges vary continuously and linearly on the numerical space, after restricting $a$ to a neighbourhood of $a_0$ we have
\[
|Z_{a,\beta}(F)-Z_{a_0,\beta}(F)| \leq\varepsilon m_{a_0,\beta}(F), \qquad 0<\varepsilon<1.
\]
Consequently,
\[
 \operatorname{Re}Z_{a,\beta}(F)
 \geq(1-\varepsilon)m_{a_0,\beta}(F)-Bd.
\]
If $F$ destabilizes at $a$, then
\[
\operatorname{Re}Z_{a,\beta}(F) \leq\frac{c^\beta(F)}d
   \operatorname{Re}Z_{a,\beta}(E) \leq\max\{0,\operatorname{Re}Z_{a,\beta}(E)\}.
\]
These inequalities bound both $m_{a_0,\beta}(F)$ and $\|v(F)\|$. Only finitely
many numerical classes can therefore destabilize locally. For each
such class, the equality of phases is linear in $a^2$. This establishes local
first-wall control using only the induction hypothesis. These bounds and the size
of the neighbourhood can be chosen uniformly on compact intervals in $\R_{>0}$, since $d$ is
fixed and $Q$ is uniformly negative on the central-charge kernels there. Thus a finite endpoint
of an interval on which $E$ is stable is an actual wall. Indeed, the uniform bound gives a finite set of
possible subobject classes near the endpoint. Continuity gives semistability there. If $E$ were still stable,
all the corresponding phase inequalities would remain strict in a neighbourhood,
contradicting that it is an endpoint.

Suppose now that $Q(E)<0$. If $E$ is strictly semistable, there is
an exact sequence
\[
  0\longrightarrow F\longrightarrow E\longrightarrow G \longrightarrow0
\]
with $F,G$ nonzero and semistable of the same phase. Both have smaller positive degree and satisfy $Q\geq0$ by induction. Their sum also satisfies
$Q\geq0$, since for every phase $\phi$,
\[
 \left\{
 v\ \middle|\
 Z_{a,\beta}(v)\in\R_{>0}e^{i\pi\phi},
 \ Q(v)\geq0
 \right\}
\]
is a convex cone. The same argument excludes a first wall as $a$ increases.
Indeed, all factors at such a wall have smaller positive
degree, and convexity would imply $Q(E)\geq0$.
Thus $E$ remains semistable for every larger $a$. The large-volume phase comparison then gives one of
three possibilities
\begin{enumerate}
 \item $E$ is a slope-semistable torsion-free module;
 \item $E$ is a torsion module;
 \item $H^{-1}(E)=F$ is slope-semistable torsion-free, and 
$H^0(E)=T$ is zero-dimensional.
\end{enumerate}
Indeed,
\[
  \lim_{a\to\infty} \frac{-\operatorname{Re}Z_{a,\beta}(E)}
{(a^2/2)\operatorname{Im}Z_{a,\beta}(E)} =-\frac{r(E)}{c^\beta(E)}.
\]
Comparing this limit in the cohomology triangle excludes two cohomology terms of positive rank, while comparison with the
slope filtration gives slope semistability of the remaining torsion-free term. In the first case, Lemma~\ref{lem:nc-slope-spread} gives
$Q(E)\geq0$. In the second, $Q(E)=c(E)^2\geq0$. In the third,
\[
Q(E)=Q(F)+2r(F)\operatorname{length}(T)\geq0.
\]
Each case contradicts $Q(E)<0$, completing the induction.

We have therefore proved
\(
Q(v(E))\geq0\) {for every semistable } $E$, and $Q|_{\ker Z_{a,\beta}}<0$. This proves the
support property on the full numerical group.
\end{proof}

\begin{remark}\label{rem:nc-numerical-kernel}
The last hypothesis is essential for the assertion of full numerical support. An arbitrary
sheaf of algebras on a surface can have numerical classes invisible to the ordinary Chern character.
\end{remark}

\subsection{Odd-dimensional intersections}
\label{subsec:odd-three-quadrics}

The numerical lattice of an odd-dimensional intersection of three quadrics makes the surface construction particularly
effective. In this case the Chern character of the underlying Clifford module retains every numerical direction. Thus the
surface charge constructed above has full numerical support, without an additional mass-Hom estimate.

\begin{theorem}\label{thm:odd-three-quadrics}
Let $m\geq1$ and let
\(
X\subset\mathbb P^{2m+4}_{\mathbb C}
\)
be a smooth complete intersection of three quadrics. Define
\[
\Ku(X)= \langle\mathcal O_X,\mathcal O_X(1),\ldots, \mathcal O_X(2m-2)\rangle^\perp.
\]
Then $\Ku(X)$ admits a numerical Bridgeland stability condition
$\widetilde\sigma_X$ with the support property on its full numerical Grothendieck
group.
\end{theorem}

\begin{proof}
Set $n=2m+1$. Adjunction gives
\(
\omega_X\simeq\mathcal O_X(2-n),
\)
so $X$ is Fano. By Lemma~\ref{lem:net-geometric-checks}, the
universal quadric bundle has smooth total space. Tensoring the decomposition of
\cite[Theorem~5.5]{KuzQuadrics2008} by $\mathcal O_X(-1)$ matches the line-bundle block with our
convention and gives an equivalence
\(
 \Psi\colon\Ku(X)\xrightarrow{\sim}
 \mathrm D^b(\mathbb P^2,\mathcal B_0).
\)
We use the canonical anti-involution of $\mathcal B_0$ to identify the right-module convention there with our left-module convention. Note that the
algebra $\mathcal B_0$ is locally free of positive rank over $\mathbb P^2$. Since $n$ is odd, weak Lefschetz and Poincar\'e duality give
\[
H^{\mathrm{even}}(X,\Q) =\bigoplus_{j=0}^n\Q h^j.
\]
The classes
\[
  \operatorname{ch}(\mathcal O_X(k)) =e^{kh} =\sum_{j=0}^n\frac{k^j}{j!}h^j, \qquad 0\leq
   k\leq n,
\]
form a basis by the Vandermonde determinant. Riemann--Roch identifies the Euler pairing with
\[
(u,w)\longmapsto\int_Xu^\vee w\,\operatorname{td}(X).
\]
This pairing is nondegenerate on $H^{\mathrm{even}}(X,\Q)$, since dualization
and multiplication by $\operatorname{td}(X)$ are invertible and the Poincar\'e
pairing is nondegenerate. Hence
\(
 \operatorname{rk}K_{\mathrm{num}}(X)=n+1.
\) Euler adjunction shows that the admissible inclusion of
$\Ku(X)$ and its adjoints preserve numerical equivalence.
The semiorthogonal decomposition therefore induces
\[
 K_{\mathrm{num}}(X)_\Q \simeq K_{\mathrm{num}}(\Ku(X))_\Q\oplus\Q^{n-2}.
\]
Consequently,
\[
\operatorname{rk}K_{\mathrm{num}}(\Ku(X)) =(n+1)-(n-2)=3,
\]
and the Clifford category has the same numerical rank. Let
$\operatorname{Forg}$ be the forgetful functor. Its exact left adjoint is
$A\mapsto\mathcal B_0\otimes A$, so
\(
\chi_{\mathbb P^2}(A,\operatorname{Forg}E) =\chi_{\mathcal B_0}(\mathcal B_0\otimes A,E).
\)
Thus numerical triviality of $E$ implies numerical triviality of $\operatorname{Forg}E$. The underlying Chern character therefore induces
\[
 v\colon K_{\mathrm{num}}(\mathbb P^2,\mathcal B_0)_\Q
 \longrightarrow\Q^3,
 \qquad
 [E]\longmapsto
 \bigl(r(E),H\cdot\operatorname{ch}_1(E),
              \operatorname{ch}_2(E)\bigr).
\]

Write $v(\mathcal B_0)=(R,C,S)$, with $R>0$. Then
\[
v(\mathcal B_0(k)) =\left(R,C+kR,S+kC+\frac{k^2R}{2}\right), \qquad k=0,1,2,
\]
and
\[
\det
 \begin{pmatrix}
R&C&S\\
 R&C+R&S+C+R/2\\
 R&C+2R&S+2C+2R
 \end{pmatrix}
 =R^3\ne0.
\]
Hence $v$ is an isomorphism. Choose $N\geq0$ such that $\mathcal B_0(N)$ is globally
generated. Proposition~\ref{prop:nc-surface-stability}, with $(a,\beta)=(1,0)$, gives a stability
condition whose support form is
\[
Q(r,c,s)=c^2-2rs+N^2r^2.
\]
Through $v$, this form gives support on the full numerical
group. Pulling the stability condition back through $\Psi$
proves the theorem.
\end{proof}

\begin{proof}[Proof of Theorem~\ref{thm:intro-three-quadrics}]
Use Theorem~\ref{thm:odd-three-quadrics} when $n$ is odd and Theorem~\ref{thm:even-three-quadrics} when $n$ is even.
\end{proof}

\section{Connected components and global coordinates}
\label{sec:stability-components}

In this section, we compare numerical and actual Grothendieck groups \cite{Vial2013} and
identify each component of stability condition spaces of odd cubics and quartic double fivefolds with
$\C\times\HH$. We then use spherical objects and twists to distinguish infinitely many components in the
two general families \cite{IlievManivel2014,Abuaf2020,SeidelThomas2001}.

\subsection{The rank-two argument}

We first record an elementary consequence of having a classical generator. The
Grothendieck group in the following lemma is the actual group, before passing to numerical
equivalence.

\begin{lemma}
\label{lem:no-real-rank-one}
Let $\D\neq0$ be an idempotent-complete triangulated
category with a classical generator, and suppose that
$K_0(\D)$ is not finitely generated. Every locally finite
stability condition whose central charge factors through a
rank-two lattice has a central charge of real rank two.
\end{lemma}

\begin{proof}
Write $\widetilde\sigma=(\cP,Z)$. Its central charge cannot have rank zero,
since a nonzero object has a nonzero semistable HN factor, and the central charge of such a
factor is nonzero. Suppose that its real image is a line. Choose $\beta\in\R$ such that this line is $\R
e^{\pi\ii\beta}$. By definition, 
\[
 \cP(\phi)\neq0\quad\Longrightarrow\quad
                         \phi\in\beta+\Z.
\]
Consequently
\[
\cP\bigl((\beta-\tfrac12,\beta+\tfrac12]\bigr) =\cP(\beta).
\]
This is the heart of a bounded $t$-structure. By local
finiteness, an interval of sufficiently small width around
$\beta$ has finite length. That interval contains no other
nonzero phase category, so $\cP(\beta)$ itself is a length
category. Let $G$ be a classical generator, and let $\mathcal S$
be the finite set of simple factors of its nonzero cohomology
objects for this bounded $t$-structure. The full subcategory of
objects whose cohomology has all its simple factors in
$\mathcal S$ is thick, by the long exact cohomology sequence
and closure under subquotients and direct summands. It contains
$G$, hence equals $\D$. Thus $\cP(\beta)$ has finitely many
simple objects, and
\[
 K_0(\D)\simeq K_0(\cP(\beta))\simeq\Z^{|\mathcal S|},
\]
a contradiction.
\end{proof}

\begin{proposition}
\label{thm:rank-two-components}
Let $\D\neq0$ be an idempotent-complete triangulated category with a classical generator, and let $v\colon K_0(\D)\twoheadrightarrow\Lambda$ be a
class map to a free lattice of rank two. Suppose that $K_0(\D)$ is not finitely generated. Every connected component $\Sigma$ of $\Stab_{\Lambda}(\D)$ is
a free $\widetilde{\mathrm{GL}}_2^+(\R)$-orbit. In particular, each connected component is contractible.
\end{proposition}

\begin{proof}
By Lemma~\ref{lem:no-real-rank-one}, every central charge in
$\Stab_\Lambda(\D)$ is a real isomorphism $\Lambda_\R\to\C$. Fix
$\widetilde\sigma=(\cP,Z)$. For $Z'$ sufficiently close to $Z$, set
\[
T=Z\circ(Z')^{-1}\in\mathrm{GL}_2^+(\R), \qquad T^{-1}Z=Z'.
\]
Let $g$ be its lift near the identity. Then $\widetilde\sigma[g]$ has charge
$Z'$ and lies in the deformation neighbourhood of
$\widetilde\sigma$. Uniqueness in the local central-charge chart
shows that the orbit contains this neighbourhood. Thus every orbit
is open. Since $\widetilde{\mathrm{GL}}_2^+(\R)$ is connected, each orbit is
connected. Its complement is a union of open orbits, so it is also closed. Hence the
orbits are precisely the connected components. If $\widetilde\sigma[g]=\widetilde\sigma$
for $g=(T,f)$, then
\[
T^{-1}Z=Z\Longrightarrow T=\id.
\]
Thus $f(\phi)=\phi+2k$ for some $k\in\Z$, and equality of slicings gives
\(
 \cP(\phi)=\cP(\phi+2k).
\)
Choose $0\ne E\in\cP(\phi)$. If $k\ne0$, phase orthogonality would force $\Hom(E,E)=0$, a contradiction. Therefore $k=0$ and $g=\id$. Let $\Sigma$
be the component containing $\widetilde\sigma$. Choose an integral basis $v_1,v_2$ of $\Lambda$ with
\(
\operatorname{Im}({Z(v_2)}/{Z(v_1)})>0.
\)
The central charges of this orientation form
\[
 \mathcal U_+
 =\left\{(z_1,z_2)\in\C^2\ \middle|\
 z_1\ne0,\ \operatorname{Im}(z_2/z_1)>0\right\}.
\]
The map
\[
 \mathcal U_+\xrightarrow{\sim}\C^*\times\HH, \qquad
  (z_1,z_2)\longmapsto(z_1,z_2/z_1)
\]
is biholomorphic. In addition, the orbit map
\[
  \widetilde{\mathrm{GL}}_2^+(\R)\longrightarrow\Sigma, \qquad g\longmapsto\widetilde\sigma[g]
\]
is a bijective local homeomorphism. Under this identification, the central-charge map
$\Sigma\to\mathcal U_+$ is the universal covering. Its local charts therefore give global
holomorphic coordinates
\begin{equation}
 u=\log Z(v_1),\qquad
 \tau=\frac{Z(v_2)}{Z(v_1)}\in\HH,
 \qquad
 Z(v_1)=e^u,\quad Z(v_2)=e^u\tau,
 \label{eq:rank-two-global-coordinates}
\end{equation}
where the logarithm is lifted after choosing its value
at one point. Hence
\(
\Sigma\simeq\C\times\HH. \) Consequently $\Sigma/\C\simeq\HH$, and
both $\Sigma$ and its quotient are contractible.
\end{proof}

\begin{remark}
\label{rem:component-versus-whole}
Proposition~\ref{thm:rank-two-components} describes every component,
but does not identify distinct components. The indexing set in
\[
\Stab_{\Lambda}(\D)\simeq \coprod_{\alpha\in\pi_0\Stab_{\Lambda}(\D)} (\C\times\HH)
\]
is part of the global topology.
\end{remark}

\subsection{Grothendieck groups of odd-dimensional examples}

The following application of \cite[Theorem~5]{Vial2013}
distinguishes the actual Grothendieck group from its
numerical quotient.

\begin{lemma}
\label{lem:odd-cohomology-kzero}
Let $X$ be a smooth projective variety over $\C$ with $H^{2j+1}(X,\Q)\neq0$ for some $j$. Suppose that
\[
\Db(X)=\langle\A,E_1,\ldots,E_r\rangle
\]
is a semiorthogonal decomposition with finitely many exceptional objects. Then $K_0(\A)\otimes\Q$ is
infinite-dimensional. In particular, $K_0(\A)$ is not finitely generated.
\end{lemma}

\begin{proof}
The rational Chern character gives an isomorphism
\[
K_0(\Db(X))_\Q \xrightarrow{\sim}\bigoplus_i\mathrm{CH}^i(X)_\Q.
\]
If this space were finite-dimensional, all rational Chow groups of $X$ would be
finite-dimensional. By \cite[Theorem~5]{Vial2013}, its
rational Chow motive would be a direct sum of Tate
motives, implying
\(
H^{\mathrm{odd}}(X,\Q)=0.
\)
This contradicts the assumption. Hence $\dim_\Q K_0(\Db(X))_\Q=\infty$. The semiorthogonal decomposition induces
\[
 K_0(\Db(X))_\Q\simeq K_0(\A)_\Q\oplus\Q^r.
\]
Therefore $\dim_\Q K_0(\A)_\Q=\infty$, proving the assertion.
\end{proof}

\begin{proposition}
\label{prop:cubic-numerical-rank}
Let $X\subset\PP^{n+1}$ be a smooth cubic hypersurface of odd dimension $n\geq3$. Then
\[
\rk\Knum(\Ku(X))=2, \qquad \dim_{\Q}\bigl(K_0(\Ku(X))\otimes\Q\bigr)=\infty.
\]
Consequently every component of $\Stab_{\mathrm{num}}(\Ku(X))$ has the global coordinates
\eqref{eq:rank-two-global-coordinates}.
\end{proposition}

\begin{proof}
Weak Lefschetz and Poincar\'e duality give
\[
H^{\mathrm{even}}(X,\Q) =\bigoplus_{j=0}^n\Q h^j \simeq\Q[h]/(h^{n+1}),
\]
where $h$ is the hyperplane class. The characters
\[
 \operatorname{ch}(\OO_X(k))
 =e^{kh}
 =\sum_{j=0}^n\frac{k^j}{j!}h^j,
 \qquad 0\leq k\leq n,
\]
form a basis by the Vandermonde determinant. Riemann--Roch identifies the Euler pairing with
\[
(u,v)\longmapsto\int_Xu^\vee v\,\operatorname{td}(X).
\]
This pairing is nondegenerate, since dualization and multiplication
by $\operatorname{td}(X)$ are invertible and the
intersection pairing is nondegenerate. Hence
\(
\rk\Knum(\Db(X))=n+1. \) Euler adjunction shows that the admissible
inclusion and projection preserve numerical equivalence. The
exceptional block has length $n-1$, so
\[
  \Knum(\Db(X))_\Q \simeq\Knum(\Ku(X))_\Q\oplus\Q^{n-1}, \qquad \rk\Knum(\Ku(X))=2.
\]

Since $n$ is odd, the middle Betti number is
\[
b_n(X)=n+1-\chi_{\mathrm{top}}(X) =\frac{2^{n+2}-2}{3}>0,
\]
where the Euler characteristic is computed from
\[
c(T_X)=\frac{(1+h)^{n+2}}{1+3h}.
\]
Lemma~\ref{lem:odd-cohomology-kzero} therefore implies that $K_0(\Ku(X))$ is not finitely generated. Finally, $\Ku(X)$ is smooth and proper and has
a classical generator. Its numerical stability manifold is nonempty by Theorem~\ref{thm:intro-hypersurfaces}. Proposition~\ref{thm:rank-two-components}
now applies.
\end{proof}

\begin{proposition}
\label{prop:double-fivefold-components}
Let $Y$ be a smooth quartic double fivefold. Then
\[
 \rk\Knum(\Ku(Y))=2, \qquad \dim_{\Q}\bigl(K_0(\Ku(Y))\otimes\Q\bigr)=\infty.
\]
Its numerical stability manifold is nonempty, and every connected component is biholomorphic
to $\C\times\HH$.
\end{proposition}

\begin{proof}
Write $Y=\{z^2+f_4=0\}\subset\PP(1,1,1,1,1,1,2)$ as in
Section~\ref{subsec:quartic-double-fivefolds}. Weighted weak Lefschetz
\cite[Section~4.2.2]{Dolgachev1982} and Poincar\'e
duality give
\[
H^{\mathrm{even}}(Y,\Q)=\bigoplus_{j=0}^5\Q h^j.
\]
The argument of Proposition~\ref{prop:cubic-numerical-rank}
therefore yields
\[
\rk\Knum(\Db(Y))=6, \qquad \rk\Knum(\Ku(Y))=6-4=2.
\]

Let
\[
  J=\C[x_1,\ldots,x_6,z]/ \bigl(\partial_{x_1}f_4,\ldots,\partial_{x_6}f_4,2z\bigr)
\]
be the weighted Jacobian ring. By the calculation in \cite[Sections~4.2.3
and~4.3]{Dolgachev1982}, we have
\[
  H^{4,1}(Y)\simeq J_{2\cdot4-(6+2)} =J_0\simeq\C.
\]
Hence $H^5(Y,\Q)\ne0$, and Lemma~\ref{lem:odd-cohomology-kzero} implies that $K_0(\Ku(Y))$ is not finitely generated.
Finally, $\Ku(Y)$ is smooth and proper and has a classical generator. Its numerical stability manifold is nonempty by
Proposition~\ref{prop:quartic-double-existence}. Proposition~\ref{thm:rank-two-components} now applies.
\end{proof}

\begin{remark}
For a cubic fivefold, the orbit constructed in \cite[Theorems~1.1
and~5.5]{PeizeLiu2025} is an entire connected component by
Proposition~\ref{prop:cubic-numerical-rank}. This conclusion also uses
the independent existence theorem of that paper. It does not
identify this component with the one obtained by our
coordinate-cut construction. The uniqueness question in
\cite[Remark~5.6]{PeizeLiu2025} concerns Serre-invariant orbits, and is
separate from the component description above. For cubic
threefolds, the distinguished component is described in
\cite[Theorem~1.1]{FanLiuMa2023}. For very general cubic fourfolds,
the full numerical stability manifold is already known
to be one $\widetilde{\mathrm{GL}}_2^+(\R)$-orbit
\cite[Corollary~3.9]{FeyzbakhshPertusi2023}.
\end{remark}

\subsection{Spherical twists and distinct components}

We use spherical objects and their twists in the sense of \cite{SeidelThomas2001}. For a
proper enhanced category with Serre functor, a $d$-spherical object $S$ satisfies
\[
\Hom(S,S[j])=
 \begin{cases}
  \C&j=0,d,\\
0&j\neq0,d,
 \end{cases}
\qquad \Sfun(S)\simeq S[d].
\]
Its twist is the autoequivalence defined by the evaluation triangle
\begin{equation}
\RHom(S,E)\otimes S\longrightarrow E \longrightarrow T_S(E).
 \label{eq:spherical-evaluation}
\end{equation}
In particular, $T_S(S)\simeq S[1-d]$.

\begin{proposition}
\label{thm:spherical-distinct-components}
Let $\D$ be a proper enhanced idempotent-complete category with
Serre functor and numerical rank two. Suppose that every connected
component of $\Stab_{\mathrm{num}}(\D)$ is a free
$\widetilde{\mathrm{GL}}_2^+(\R)$-orbit. If $S\in\D$ is
$d$-spherical for $d\geq2$, then
\(
 T_S^m(\Sigma)\neq\Sigma
\)
for every connected component $\Sigma$ and $m\ne0$. In particular, if the numerical
stability manifold is nonempty, it has infinitely many connected components.
\end{proposition}

\begin{proof}
It suffices to consider $m>0$, since a negative power preserves a component if and only if its inverse
does. Suppose $T_S^m(\Sigma)=\Sigma$, and fix $\widetilde\sigma=(\cP,Z)\in\Sigma$. Since
$\Sigma$ is a $\widetilde{\mathrm{GL}}_2^+(\R)$-orbit, there is $g=(T,f)$ such that
\(
T_S^m\widetilde\sigma=\widetilde\sigma[g].
\)
Thus $T_S^m$ takes semistable objects of phase $\phi$ to semistable objects of phase $f(\phi)$, where $f$
{ is strictly increasing}, and $f(\phi+1)=f(\phi)+1$. Let $\beta_1>\cdots>\beta_l$ be the HN
phases of $S$. Applying $T_S^m$ to its HN filtration and using
\(
T_S^m(S)\simeq S[-m(d-1)]
\)
gives, by uniqueness,
\begin{equation}
f(\beta_j)=\beta_j-c, \qquad c=m(d-1)\geq1.
 \label{eq:spherical-phase-displacement}
\end{equation}
Choose one phase $\beta=\beta_j$. For $\beta\leq\phi<\beta+1$,
\[
f(\phi)<f(\beta+1)=\beta-c+1\leq\beta\leq\phi.
\]
Periodicity therefore implies
\begin{equation}
 f(\phi)<\phi\text{ for }\phi\in\R.
 \label{eq:spherical-all-phases-decrease}
\end{equation}

For every $\widetilde\sigma$-semistable object $E$, phase orthogonality now gives
\(
 \Hom(E,T_S^mE)=0.
\) Consider the composite of the natural morphisms in
\eqref{eq:spherical-evaluation},
\[
E\longrightarrow T_SE\longrightarrow\cdots \longrightarrow T_S^mE.
\]
Each successive cone belongs to $\thick(S)$, so the octahedral axiom gives the
same conclusion for the cone of the composite. Since the composite vanishes, we obtain
\(
T_S^mE\oplus E[1]\in\thick(S).
\)
Closure under direct summands and shifts implies $E\in\thick(S)$. Finite HN filtrations then give
\(
\D=\thick(S). \) Set $A=\RHom(S,S)$. Differential graded Morita theory identifies $\D$ with $\Perf(A)$,
and sphericity gives
\[
H^j(A)=
 \begin{cases}
\C,&j=0,d,\\
  0,&\text{otherwise}.
 \end{cases}
\]
By \cite[Theorem~8.1]{KellerNicolas2013}, $\Perf(A)$ has a bounded $t$-structure with a
length heart whose simple objects correspond to the indecomposable summands of $A$. Since $H^0(A)=\C$,
there is exactly one simple object. Hence
\[
  K_0(\D)\simeq K_0(\Perf(A))\simeq\Z, \qquad \rk\Knum(\D)\leq1,
\]
contrary to $\rk\Knum(\D)=2$. Thus no nonzero power of $T_S$ preserves a
component. In particular,
\[
  T_S^i(\Sigma)=T_S^j(\Sigma)\Longrightarrow i=j,
\]
so $j\mapsto T_S^j(\Sigma)$ defines an injection
\(
\Z\hookrightarrow\pi_0\Stab_{\mathrm{num}}(\D).
\)
\end{proof}

The required spherical objects in the following corollary are geometric constructions from
\cite{IlievManivel2014,Abuaf2020}.

\begin{corollary}
\label{cor:sevenfold-double-fivefold-infinite}
For a general cubic sevenfold $X$ and a general quartic double fivefold $Y$,
the manifolds $\Stab_{\mathrm{num}}(\Ku(X))$ and $\Stab_{\mathrm{num}}(\Ku(Y))$ have
infinitely many connected components. Every component is biholomorphic to
$\C\times\HH$, and every nonzero power of a suitable spherical twist
moves every component.
\end{corollary}

\begin{proof}
For a general cubic sevenfold $X$, the Cartan representation gives a rank-nine vector bundle
$\mathcal E_X$ such that
\(
 \mathcal E_X(-1),\ \mathcal E_X(-2)\in\Ku(X)
\)
are $3$-spherical by \cite[Proposition~5.2, Theorem~3.4 and
Corollary~5.3]{IlievManivel2014}.
Proposition~\ref{prop:cubic-numerical-rank} then gives
the rank-two component description, so
Proposition~\ref{thm:spherical-distinct-components} applies
to either object.

For a general quartic double fivefold $Y$, the rank-six bundles
\(
 \widetilde{\mathcal E}(-1),\
 \widetilde{\mathcal E}(-2)\in\Ku(Y)
\)
are $3$-spherical by \cite[Theorem~3.8]{Abuaf2020}.
The conclusion follows in the same way from
Proposition~\ref{prop:double-fivefold-components}
and Proposition~\ref{thm:spherical-distinct-components}.
\end{proof}

\begin{proof}[Proof of Theorem~\ref{thm:intro-cubic-components}]
Propositions~\ref{prop:cubic-numerical-rank} and~\ref{prop:double-fivefold-components} identify each connected component and its $\C$-quotient.
Their disjoint unions give the asserted global descriptions. For the two general families,
Corollary~\ref{cor:sevenfold-double-fivefold-infinite} gives infinitely many components and shows that no nonzero power of the
spherical twists preserves a component.
\end{proof}

\begin{remark}
The fractional Calabi--Yau identity for a degree-$d$ hypersurface is
\(
 \Sfun^d\simeq[(n+2)(d-2)].
\)
If a nonzero object is $s$-spherical, this identity forces $ds=(n+2)(d-2)$. A nonzero object
in a proper category cannot be isomorphic to a nonzero shift of itself,
since iteration would give nonzero morphisms in arbitrarily
large degrees. In the cubic case, $3s=n+2$. In particular, cubic threefolds
and fivefolds have no spherical objects. This observation
does not determine their connectedness. For ordinary quartics,
sextics, and the intersections of three quadrics
considered above, numerical rank is in general larger
than two, and the rank-two argument does not supply
their full component description.
\end{remark}

\begin{remark}
\label{rem:moduli-consequences}
The mass-Hom bounds established above also give consequences for moduli spaces.
For the numerical stability conditions $\widetilde\sigma$ satisfying these
bounds, the stack of semistable objects of fixed numerical class $v$ and phase $\phi$ admits
a proper good moduli space $M_{\widetilde\sigma}^{\mathrm{ss}}(v,\phi)$ by \cite[Example~2.29, Remark~2.36
and Theorem~2.35]{HLR2025}. The categories here are smooth and proper, and numerical
stability over $\C$ is locally constant in the sense of that theorem.

In the rank-two examples considered in this section, stable and semistable objects are
unchanged throughout each connected component, up to reindexing their phases. More precisely, for
$g=(T,f)\in\widetilde{\mathrm{GL}}_2^+(\R)$ there are natural isomorphisms
\(
M_{\widetilde\sigma[g]}^{\mathrm{ss}} (v,f^{-1}(\phi)) \simeq
M_{\widetilde\sigma}^{\mathrm{ss}}(v,\phi).
\)
Moreover, every semistable object of primitive numerical class is stable. Indeed, the central charge is a real isomorphism, so the numerical
classes of its Jordan--H\"older factors are positive integral multiples of that primitive class.

In the examples of Corollary~\ref{cor:sevenfold-double-fivefold-infinite}, spherical twists
similarly identify moduli spaces associated with distinct components when both the stability
condition and the numerical class are transported. These observations concern properness and invariance of
moduli spaces. Nonemptiness for every nonzero numerical class and the Lagrangian restriction
constructions of \cite[Theorems~7.5 and~7.18 and Corollary~7.20]{PeizeLiu2025}
require additional geometric arguments.
\end{remark}

\end{document}